\documentclass[11pt,reqno]{amsart}

\usepackage[utf8]{inputenc}
\usepackage[margin=1.25in]{geometry}
\usepackage{hyperref}
\usepackage{appendix}
\usepackage{amsfonts}
\usepackage{amssymb}
\usepackage{stmaryrd} 
\usepackage{amsmath}
\usepackage{amsthm}
\usepackage{mathrsfs}

\theoremstyle{plain}
\newtheorem{theorem}{Theorem}[section]
\newtheorem{claim}[theorem]{Claim}
\newtheorem{corollary}[theorem]{Corollary}
\newtheorem{lemma}[theorem]{Lemma}
\newtheorem{proposition}[theorem]{Proposition}
\newtheorem{question}[theorem]{Question}

\theoremstyle{remark}

\newtheorem{remark}[theorem]{Remark}

\usepackage{biblatex} 
\numberwithin{equation}{section}

\title[Smoothness density McKean-Vlasov SDEs]{Smoothness of the density for McKean-Vlasov SDEs with measurable kernel}
\author{Yi Han}
\address{Department of Pure Mathematics and Mathematical Statistics, University of Cambridge.
}
\email{yh482@cam.ac.uk}
\thanks{Supported by EPSRC grant EP/W524141/1.}

\begin{document}

\begin{abstract}
Consider the McKean-Vlasov SDE $$
    dX_t=\langle b(X_t-\cdot),\mu_t\rangle dt+dW_t,\quad \mu_t=\operatorname{Law}(X_t),
$$
where $W$ is the $n$-dimensional Brownian motion and $b:\mathbb{R}^d\to\mathbb{R}^d$ is a measurable function. First assuming $b\in L^\infty$, we prove that the law $\mu_t$ of $X_t$ has a density $p_t$ with respect to the Lebesgue measure, which is continuously differentiable with gradient being $\gamma$-Hölder continuous for each $\gamma\in(0,1)$. Assume further that $b\in \mathcal{C}_b^{1}$, we prove that the density $p_t$ is infinitely differentiable. In the regularization by noise perspective, this shows McKean-Vlasov SDEs tend to have a smoother density function than SDEs without density dependence, under the same regularity assumption of the coefficients. We observe similar phenomenon for singular interaction kernels satisfying Krylov's integrability condition, for distributional kernels $b\in B_{\infty,\infty}^\alpha$, $\alpha\in(-1,0)$, and for processes driven by an $\alpha$-stable noise for $\alpha\in(1,2)$.
\end{abstract}

\maketitle

\section{Introduction}

In this paper, we are interested in regularity of probability density functions for McKean-Vlasov SDEs 
\begin{equation}
 \label{SDEMCKEAN}   
dX_t=\langle b(X_t-\cdot),\mu_t\rangle dt+dW_t,\quad \mu_t=\operatorname{Law}(X_t)\end{equation} with a nonsmooth interaction kernel $b:\mathbb{R}^d\to\mathbb{R}^d$. 

Smoothness of density functions for diffusion processes is a central topic in stochastic analysis, with Malliavin calculus being one of the key techniques. Consider a diffusion process 
$$dX_t=V(X_t)dt+dW_t,$$
then as a special case of the celebrated Hörmander's theorem \cite{hormander1967hypoelliptic} \cite{malliavin1978stochastic} \cite{norris1986simplified}, a smooth vector field $V$ and a nondegenerate Brownian noise imply that the diffusion process has a smooth density with respect to the Lebesgue measure. On the other hand, if the vector field $V$ is nonsmooth, then we in general do not have a smooth density function. More details in this direction will be discussed in Section \ref{sec1.2}.

Malliavin calculus, as the central toolbox in proving regularity of densities, works when the vector field $V$ is at least once continuously differentiable. In recent years, new techniques have been invented to prove regularity of density functions when $V$ is not differentiable (see \cite{debussche2014existence}, \cite{debussche2013existence}, \cite{romito2018simple} and \cite{bally2017convergence}). These techniques take full advantage of the driving Brownian motion $W$, and use a linearisation of the diffusion process $X_t$ at time $t-\epsilon$ to approximate. However, the regularity of densities in these cases are rather weak, as most of the times they are not proven to be continuously differentiable. More details will be discussed in section \ref{sec1.1}.

In this paper, we study density dependent SDEs \eqref{SDEMCKEAN} taking into account both the regularizing properties of the Brownian motion $W$ and smoothing properties of the convolution drift term $\langle b(x-\cdot),\mu_t\rangle.$

\subsection{Regularization by noise and smoothing via convolution}
To begin the story, we give a brief account of the regularizing effects of Brownian motion. Consider a vector field $V:\mathbb{R}^d\to\mathbb{R}^d$ which is merely bounded measurable. Then the ODE $$dX_t=V(X_t)dt$$ is in general not well-posed, as the Cauchy-Lipschitz theorem fails. However, imposing a Brownian noise, the SDE
\begin{equation}\label{sdetype}dX_t=V(X_t)dt+dW_t,\end{equation} has a unique weak solution thanks to Girsanov's theorem. Strong well-posedness also holds, see for example \cite{zvonkin1974transformation}, \cite{menoukeu2013variational} and \cite{menoukeu2019flows}.

There has been a growing interest in McKean-Vlasov SDEs of the form \eqref{SDEMCKEAN}. A review for its applications in physics and modelling can be found in \cite{chaintron2021propagation}, and applications in stochastic control problem can be found in \cite{carmona2018probabilistic}. In recent years, much attention has been paid to the case where the interaction kernel $b$ is nonsmooth \cite{jabin2016mean}. It is by now well-understood that the regularity assumption on $b$ for \eqref{SDEMCKEAN} to be well-posed, is at least as general as the regularity assumptions on the vector field $V$ for the SDE\eqref{sdetype} to be well-posed. This includes in particular kernels $b$ of linear growth, or satisfying the integrability condition named after Krylov. We refer to \cite{han2022solving} for a quick approach in proving well-posedness of \eqref{SDEMCKEAN} under general regularity assumptions, via a fixed point problem in relative entropy.

The smoothing property of convolution by distribution $\langle b(x-\cdot,\mu\rangle$ is brought to light only very recently.
In a very recent work \cite{de2022multidimensional}, it is proved that McKean-Vlasov SDEs \eqref{SDEMCKEAN} could be solved with $b$ belonging to a much wider regularity class, such that with $V$ in the same class, the SDE \eqref{sdetype} is not yet known to be well-posed. This involves distributional kernels $b\in\mathcal{C}_b^\alpha$ for some $\alpha\in(-1,-\frac{1}{2}].$ The reason is that convolution by density $\langle b(x-\cdot,\mu_t\rangle$ tends to make the drift smoother if $\mu_t$ itself has some regularity, so the drift $b(x-\cdot,\mu_t\rangle$ has better regularity than what is assumed on $b$.

In this paper, we take a different perspective towards this regularizing by convolution effect, by proving that the probability density function of \eqref{SDEMCKEAN}, and the drift $\langle b(t,X_t-\cdot),\mu\rangle$ of \eqref{SDEMCKEAN} as well, can be a lot smoother than what one may expect. A loosely related work in this direction is \cite{crisan2018smoothing}, where coefficients of the McKean-Vlasov SDE may depend on $\mu_t$ non-linearly, and they proved that the density has certain smoothness properties when coefficients are regular enough. Via a bootstrap argument, we prove the density is a lot smoother than the kernel $b$, in the  case of linear dependence $\langle b(x-\cdot),\mu_t\rangle$. Another notable work is \cite{de2021backward}, where the authors studied transition densities of a backward Kolmogorov PDE with Hölder coefficients and obtained some propagation of chaos estimates. We study the smoothing properties from a different direction by focusing on improved regularity of $\mu_t$, see Corollary \ref{corollary1.3}.

\subsection{Bounded measurable kernels}\label{sec1.1} For the purpose of comparison, we review possible behaviours of the density functions of SDE \eqref{sdetype}, only assuming $V\in L^\infty$. It is proved in \cite{romito2018simple} that a density function for $X_t$ with respect to Lebesgue measure always exists, and the density lies in some Besov spaces. In Section \ref{section2}, we will elaborate on techniques developed in \cite{romito2018simple}. To the other extreme, consider the one-dimensional SDE
\begin{equation}
    \label{badsde}
X_t=-\int_0^t \operatorname{}{sign}(X_s)ds+W_t,\quad t\geq 0.\end{equation}
For each $t>0$, \eqref{badsde} has a density function which is Lipschitz continuous on $\mathbb{R}^d$, but has no more regularity beyond that. The computational details are quoted in Appendix \ref{appendixB}. 

For McKean-Vlasov SDEs of the form \eqref{SDEMCKEAN}, we prove the following theorem:

\begin{theorem}\label{theorem1}
Fix an interaction kernel $b(t,x):[0,\infty)\times\mathbb{R}^d\to\mathbb{R}^d$ such that $b(t,x)\in L^\infty([0,\infty)\times\mathbb{R}^d)$. Consider the McKean-Vlasov SDE
\begin{equation}\label{1.43}
    dX_t=\langle b(t,X_t-\cdot),\mu_t\rangle dt+dW_t,\quad\operatorname{Law}(X_t)=\mu_t,
\end{equation} with initial distribution $\mu_0\in\mathcal{P}(\mathbb{R}^d)$, the space of Borel probability measures on $\mathbb{R}^d$.

Then at each time $t>0$, $X_t$ has a density $p_t$ with respect to the Lebesgue measure that lies in both $B_{1,\infty}^{2+\gamma}(\mathbb{R}^d)$ and $\mathcal{C}_b^{1+\gamma}(\mathbb{R}^d)$ for every $\gamma\in(0,1).$
\end{theorem}

The well-posedness of \eqref{1.43} can be found in \cite{lacker2018strong}. This theorem shows $p_t$ is continuously differentiable with a Hölder continuous gradient, which is apparently more regular than that of the SDE \eqref{badsde}. Moreover, the drift of \eqref{1.43} is $\mathcal{C}^2$ at each time $t>0$ as opposed to merely $L^\infty$, which might be a bit surprising at first glance.

\subsection{Bootstrapping the regularity}
\label{sec1.2} Now we assume the kernel $b$ has some regularity, we show that the regularity can be propagated and give rise to a smooth density.

\begin{theorem}\label{theorem2}
Fix an interaction kernel $b(t,x):[0,\infty)\times\mathbb{R}^d\to\mathbb{R}^d$ such that we have $b(t.\cdot)\in \mathcal{C}_b^{1}(\mathbb{R}^d)$ (the space of differentiable functions with bounded derivative, with the norm $\|f\|_{\mathcal{C}_b^1}:=\|f\|_{L^\infty}+\|\nabla f\|_{L^\infty}$) for each $t>0$, satisfying $\|b(t,\cdot)\|_{\mathcal{C}_b^{1}}\in L^\infty([0,\infty))$.
Consider the McKean-Vlasov SDE
\begin{equation}\label{densitydependence}
    dX_t=\langle b(t,X_t-\cdot),\mu_t\rangle dt+dW_t,\quad\operatorname{Law}(X_t)=\mu_t
\end{equation}
with initial law $\mu_0\in\mathcal{P}(\mathbb{R}^d)$ and $W$ is the $d$-dimensional Brownian motion. Then at each time $t>0$, $X_t$ has a \textbf{smooth} density $p_t$ with respect to the Lebesgue measure.
\end{theorem}

We make a comparison with the SDE case \eqref{sdetype}. It is proved in \cite{banos2016malliavin} that if $V\in\mathcal{C}^k$ for some $k\in\mathbb{N}^+$, then the density at time $t$ of $dX_t=V(X_t)dt+dW_t$ lies in $\mathcal{C}^{k-1+\alpha}$ for each $\alpha\in(0,1)$, i.e. the density is $k-1$-times differentiable with $k-1$-th gradient $\alpha$-Hölder continuous, for all $\alpha\in(0,1)$. It is unclear if this regularity is sharp. 

In the setting of Theorem  \ref{theorem2}, we also have that the drift of \eqref{densitydependence} is indeed $\mathcal{C}^\infty$ at each time $t>0$.

\subsection{Sketch of proof and open questions}

The proof of both theorems rely on estimating the densities in some Bssov space norms. Definitions and some useful properties of Besov spaces are given in Section \ref{section2.1}.

The proof of Theorem \ref{theorem1} occupies Section \ref{section2} and \ref{section4}. Since $b\in L^\infty,$ the result in \cite{romito2018simple} implies $p_t\in B_{1,\infty}^\alpha$ for some $\alpha\in(0,1)$, then by convolution the drift is $\alpha$-Hölder. Then apply the technique twice, one shows $p_t$ lies in $B_{1,\infty}^{1+\beta}$ and then in $B_{1,\infty}^{2+\gamma}$ for any $\beta,\gamma\in(0,1).$ Then the conclusion of Theorem \ref{theorem1} is justified in dimension one via Sobolev embedding. In higher dimensions, we use Malliavin calculus once and then resort to a technique based on Riesz transform, which is developed in \cite{bally2011riesz}, to improve the regularity of our estimate. This last part is done in Chapter \ref{section4}. 

The proof of Theorem \ref{theorem2} occupies Section \ref{section3}. Since $b$ is continuously differentiable, we can apply Malliavin calculus once and deduce that the density function $p_t$ lies in $B_{1,1}^1$. Then the drift $\langle b(t,x-\cdot),\mu_t\rangle$ lies in the Sobolev space $W^{2,\infty}$ and we can apply Malliavin calculus twice, showing that $p_t\in B_{1,1}^2$. Iterating the above procedure shows $p_t$ is smooth. In the proof we will use a nonstandard form of the smoothing lemma, see Lemma \ref{lemma3.80}, which shows the density lies in the Besov space $B_{1,1}^s$ with the optimal value $s$ one can expect. 

We give the following corollary, which can be seen as an intermediate case of the results of Theorem \ref{theorem1} and \ref{theorem2}.

\begin{corollary}\label{corollary1.3}
Consider the McKean-Vlasov SDE \eqref{1.43}, assuming $b(t,\cdot)\in \mathcal{C}_b^\alpha$ for some $\alpha>0$, and that $\sup_{t\geq 0}\|b(t,\cdot)\|_{\mathcal{C}_b^\alpha}<\infty$. Then for each $t>0$, $\mu_t$ has a density $p_t$ which is in $\mathcal{C}_b^{2+\gamma}(\mathbb{R}^d)$ for each $\gamma\in(0,1)$.
\end{corollary}
The proof of this Corollary is very similar to that of Theorem \ref{theorem1} and is given in Section \ref{section4.3}. 

In the proof of both theorems, before using the standard techniques one should notice that since the norms of $p_t$ are likely to explode as $t\to 0$, the smoothness of the drift can be controlled only after a positive time. To solve this problem, we use Markov property of the SDE, conditioning on a smaller time $t_0<t$, and use probability estimates with the SDE starting at $t_0$. Therefore it is natural to raise the following question:
\begin{question}
Can we prove a similar result if the process is not Markovian? The possible cases are: interactions $b$ that are path dependent (see for example \cite{lacker2018strong}), or replace the Brownian motion $W$ by a fractional Brownian driving noise (see \cite{galeati2022distribution} or \cite{han2022solving}).
\end{question}

Both the conclusion in Theorem \ref{theorem1}, i.e. $p_t\in \mathcal{C}_b^{1+\alpha}$, and the assumption in Theorem \ref{theorem2}, i.e. $b\in \mathcal{C}_b^{1}$, are technical and are tailored to the different proof techniques of two theorems, see Remark \ref{remark1} and Remark \ref{remark2}. Thus we raise the following question:

\begin{question}
In the setting of Theorem \ref{theorem1}, i.e. if we only assume $b\in L^\infty$, is it possible to prove that the density $p_t$ is smooth?
\end{question}
A tentative solution to this question is outlined in Section \ref{section2.5}, yet the general picture is far from complete. We believe that it is possible to improve the results in Theorem \ref{theorem1} and Corollary \ref{corollary1.3} to a notable degree, but each improvement could be difficult.

From our perspective, the reason is as follows: in the proof of Theorem \ref{theorem1}, we use the observation that for a SDE with additive noise, on short time scales the process looks like the noise, which is a rather vague approximation with low precision. In the proof of Theorem \ref{theorem2}, we use the fact that we can integrate by parts many times in the Malliavin sense, which is much stronger than the previous approximation and allows us to pass to smoothness without much difficulty. We suspect that even if we can prove the density $p_t$ is smooth in the setting of Theorem \ref{theorem1}, the numerical constant would be intractable and the scale of computation could be forbidding.

\subsection{More Examples}

Then in Section \ref{section5} we explore more examples where our theme extends, i.e. where the density of McKean-Vlasov SDE can be smoother. In Section \ref{section5.1} and \ref{section5.20} we prove the following for singular interactions $b$, generalizing Theorem \ref{theorem1}:
\begin{theorem}\label{theorem3}
Assume that the kernel $b(t,x)\in L^q([0,T];L^p(\mathbb{R}^d))$ with $\frac{d}{p}+\frac{2}{q}<1$, and consider the McKean-Vlasov SDE \eqref{1.43}. Then $\mu_t$ has a density $p_t$ with respect to Lebesgue measure, $p_t\in B_{1,\infty}^{\gamma}$ for any $\gamma\in(0,3-\frac{4}{q})$. In the special case $q=\infty$ we prove $p_t\in B_{1,\infty}^{2+\gamma}(\mathbb{R}^d)$ and $p_t\in\mathcal{C}_b^{1+\gamma}(\mathbb{R}^d)$ for every $\gamma\in(0,1).$
\end{theorem}

By convolution, the drift $\langle b(t,x-\cdot),\mu_t\rangle\in B_{p,\infty}^{\gamma}\hookrightarrow B_{\infty,\infty}^{\gamma-\frac{d}{p}}$. Since we can choose $\gamma$ so that $\gamma-\frac{d}{p}>1$, we conclude that at each time $t>0$, the drift of \eqref{1.43} is Lipschitz continuous, though the Lipschitz constant tends to blow up as $t\to 0$.

In Section \ref{section5.30} we consider drifts given by a distribution in $\mathcal{C}_b^\alpha$, $\alpha\in(-1,0)$. We utilize the (very recent) well-posedness result in \cite{de2022multidimensional} and our results in Section \ref{section2} to obtain the following, which is indeed the same regularity result as Theorem \ref{theorem1}.

\begin{theorem}\label{theorem1.7}
Assume that $b(t,x)\in L^\infty([0,T];B_{\infty,\infty}^\alpha(\mathbb{R}^d))$ for $\alpha>-\frac{1}{2}$, and consider the McKean-Vlasov SDE \eqref{1.43}. Then the probability density function $p_t$ of $\mu_t$ lies in $B_{1,\infty}^{2+\gamma}(\mathbb{R}^d)$ and $\mathcal{C}_b^{1+\gamma}(\mathbb{R}^d)$ for each $\gamma\in(0,1)$. 
\end{theorem}

The case of $b(t,\cdot)\in B_{\infty,\infty}^\alpha$, $\alpha\in(-1,-\frac{1}{2}]$ is even more intriguing. Our strategy gives a hint as to why the corresponding SDE $dX_t=V(X_t)dt+dW_t$
with $V$ in this regularity class may not be well-posed. See \cite{flandoli2017multidimensional} or the end of Section \ref{section5.30} for more explanations.

In Section \ref{section5.2} we also make the observation that results close to Theorem \ref{theorem1} and \ref{theorem2} can be established for
McKean-Vlasov SDEs driven by an $\alpha$-stable process, for $\alpha\in(1,2)$. 

We summarize the results of this paper in the following table: for the McKean-Vlasov SDE \eqref{1.43}, without special mention we assume the associated norms of $b(t,\cdot)$ are uniform on $t\in[0,T]$, then we have:

\begin{center}
\begin{tabular}{|c|c|}%
\hline  
Regularity of the interaction kernel $b$ & Regularity of the probability density function $p_t$
\\
\hline  
$b\in\mathcal{C}_b^{1}$ & Smooth.\\
\hline 
$b\in\mathcal{C}_b^{\alpha}$, $\alpha>0$& $ \mathcal{C}_b^{2+\gamma}\cap B_{1,\infty}^{2+\gamma}$ for any $\gamma\in (0,1)$
.\\
\hline
$b\in L^\infty$, or $b\in L^p,p>d$
&$\mathcal{C}_b^{1+\gamma}\cap B_{1,\infty}^{2+\gamma}$ for any $\gamma\in(0,1)$. \\
\hline 

$b\in \mathcal{C}_b^\alpha$,  $\alpha> -\frac{1}{2}$
 & $\mathcal{C}_b^{1+\gamma}\cap B_{1,\infty}^{2+\gamma}$ for any $\gamma\in(0,1)$.\\ 
\hline
$b\in L^q([0,T];L^p)$, $\frac{d}{p}+\frac{2}{q}<1$ & $ B_{1,\infty}^\gamma$ for any $\gamma\in(0,3-\frac{4}{q}).$
\\
\hline
$b\in\mathcal{C}_b^\alpha$,  $\alpha>-1$ 
& $ B_{1,\infty}^\gamma$ for any $\gamma\in(0,1)$.\\
\hline
\end{tabular}
\end{center}

\begin{remark}
\begin{enumerate}
    \item  The case $\alpha\in (-1,-\frac{1}{2}]$ is addressed in \cite{de2022multidimensional}. This paper shows that, as long as $\alpha>-\frac{1}{2}$, the regularity of the drift is the same no matter $\alpha\in (-\frac{1}{2},0)$ or $\alpha=0$. And in all these cases except the $\alpha<-\frac{1}{2}$ one, the drift of the McKean-Vlasov SDE \eqref{1.43} is Lipschitz continuous in space at any $t>0$, so we can regard it as a sufficiently regular function although its Lipschitz constant my explode as $t\to 0$.
    \item For $L^q-L^p$ kernels, the temporal index $q$ has more impact on the regularity than the spatial index $p$. This time dependent case may deserve further investigation.
    \item By proving $p_t\in B_{1,\infty}^s$ for $s\in\mathbb{R}_+$, we actually prove the stronger version $p_t\in B_{1,1}^s$, and use the continuous embedding $B_{1,1}^s\hookrightarrow B_{1,\infty}^s$.
    \item It is possible, though rather lengthy, to make explicit the dependence of the norms of $p_t$ with respect to $t$. We choose not to do this as we are more interested in the regularity issue.
    \item In the setting of Theorem \ref{theorem2}, it is natural to pose the following question: if the initial distribution $\mu_0$ has a smooth density $p_0$, then do we have a convergence of the density function $p_t$ to $p_0$ as $t\to 0$, in the norm of continuous functions $\mathcal{C}_b^k$, for some $k\in\mathbb{N}^+$ large? We are not aware of a solution to this question.
\end{enumerate}
\end{remark}

\section{Low regularity regime: Stochastic Taylor Expansion}\label{section2}

In this section we work in the case $b\in L^\infty.$ Then regardless of the measure $\mu$, $x\mapsto \langle b(t,x-\cdot),\mu\rangle\in L^\infty.$ Thus we can regard the McKean-Vlasov SDE \eqref{SDEMCKEAN} as a SDE \eqref{sdetype} with $V\in L^\infty,$ and improve regularity from this perspective. 

In Section \ref{section2.1} we review basic notions of Besov spaces that will be used throughout the paper. Then in section \ref{section2.2} we establish regularity of density in $B_{1,\infty}^\alpha$; in section \ref{section2.3} we establish regularity in $B_{1,\infty}^{1+\alpha}$, and finally in Section \ref{section2.4} we establish regularity in  $B_{1,\infty}^{2+\alpha}.$ Finally, in Section \ref{section2.5} we outline a path that could lead to show the density is smooth.

\subsection{On Besov spaces}\label{section2.1}

We first recall the notion of Besov spaces. A comprehensive treatment of this topic can be found for example in \cite{sawano2018theory} or \cite{bahouri2011fourier}, and we frequently refer to Appendix A.1 of \cite{romito2018simple}. We start with the Littlewood-Paley decomposition: fix a dyadic decomposition $(\varphi_n)_{n\geq 0}$ on $\mathbb{R}^d$ as follows: choose $\psi$ and $\varphi$ in $\mathcal{C}_c^\infty$ satisfying $$\operatorname{supp}(\varphi)\subset B(4)\setminus B(1),\quad \operatorname{supp}(\psi)\subset B(4),\quad  \varphi_j:=\varphi(2^{-j}\cdot),$$ satisfying $\psi+\sum_{j=1}^\infty \varphi_j=1.$ Set $\varphi_0=\psi.$ ( See for example Definition 2.7 of \cite{sawano2018theory} for more options.) Then for any tempered distribution $f$ on $\mathbb{R}^d$ the factor $f_n=\mathcal{F}^{-1}(\varphi_n\hat{f})$ is a Schwartz function with
$f=\sum_n f_n$, where $\hat{f}$ is the Fourier transform of $f$ and $\mathcal{F}^{-1}$ is the inverse Fourier transform. Then for $p\geq 1, q\geq 1$ and $s\in\mathbb{R}$, define 
$$\|f\|_{B_{p,q}^s}:=\|(2^{ns}\|f_n\|_{L^p})_{n\geq 0}\|_{\ell^q}.$$
The Besov space $B_{p,q}^s(\mathbb{R}^d)$ is then defined as the closure of the Schwartz space with respect to the $\|\cdot\|_{B_{p,q}^s}$ norm.

For a vector-valued function $V=(V^1,\cdots,V^d):\mathbb{R}^d\to\mathbb{R}^d$, saying $V\in B_{p,q}^s$ means that each component $V^i\in B_{p,q}^s$, for $i=1,\cdots,d$.

We recall the embedding of Lebesgue spaces into Besov spaces \cite{sawano2018theory}: for any $\ell\in[1,\infty],$
\begin{equation}
    \label{functionrelations}
B_{\ell,1}^0\hookrightarrow L^\ell\hookrightarrow B_{\ell,\infty}^0.\end{equation}

The following convolution inequality will be frequently used: for some $\gamma\in\mathbb{R}$, $\ell\geq 1$, 
\begin{equation}
    \label{convolutioninequ}
\|f*g\|_{B_{\ell,\infty}^\gamma} \leq C_d \|f\|_{B_{\ell_1,\infty}^{\gamma-\delta}}\|g\|_{B_{\ell_2,\infty}^\delta}\end{equation}
for any $\delta\in\mathbb{R}$ and any $1+1/\ell=1/\ell_1+1/\ell_2$, where $C_d$ is a constant that depends only on $d.$ We have a more general form which is also quite useful:
\begin{equation}\label{convolutioneq2}
    \|f*g\|_{B_{\ell,m}^\gamma} \leq C_d \|f\|_{B_{\ell_1,m_1}^{\gamma-\delta}}\|g\|_{B_{\ell_2,m_2}^\delta}
\end{equation}
for $m_1$, $m_2$ satisfying $m_1^{-1}\geq (m^{-1}-m_2^{-1})\vee 0.$
The proof can be found in \cite{burenkov1989estimates}, Theorem 3.

The spaces $B_{\infty,\infty}^s,s\in\mathbb{R}^+\setminus\mathbb{Z}$ are of particular interest, as they are equivalent to the Hölder space $C_b^s(\mathbb{R}^d):=C_b^{\lfloor s\rfloor,s-\lfloor s\rfloor}(\mathbb{R}^d)$ of bounded measurable functions that are $\lfloor s\rfloor$-times differentiable, whose $\lfloor s\rfloor$-th derivatives are $s-\lfloor s\rfloor$-Hölder continuous on $\mathbb{R}^d$.

As suggested in \cite{romito2018simple}, we consider an alternative norm for the Besov space which are more useful for our purpose. For $h\in\mathbb{R}$ define the operators $$\Delta_h^1f(x)=f(x+h)-f(x),$$
$$\Delta_h^nf(x)=\Delta h^1(\Delta h^{n-1}f)(x)=\sum_{j=0}^n (-1)^{n-j}C_n^jf(x+jh).$$ For $m>s$ set the seminorm 
\begin{equation}\label{norms1}
    [f]_{B_{p,q}^s}:=\|h\mapsto \frac{\|\Delta_h^mf\|_{L^p}}{|h|^s}\|_{L^q(B_1(0);\frac{dh}{|h|^d})},
\end{equation}
and for $s>0$, $1\leq p,q\leq\infty$ define the following norms $$\|f\|_{L^p}+[f]_{B_{p,q}^s}.$$ These norms are equivalent to the standard $B_{p,q}^s$ norms (see \cite{triebeltheory}). In the setting of $B_{\infty,\infty}^s$ spaces with $s\in\mathbb{R}^+\setminus\mathbb{Z}$, we see that the Hölder norm $\|f\|_{\mathcal{C}_b^s}$ is equivalent to the norm 
\begin{equation}\label{supernorm}\|f\|_{L^\infty}+\sup_{|h|\leq 1, h\neq 0}\frac{\|\Delta_h^m f\|_{L^\infty}}{|h|^s},\end{equation}
whenever $m>s$.

For $s\in\mathbb{R}^+$ let $\mathscr{C}_b^s$ denote the closure of the Schwartz space with respect to the norm \eqref{supernorm}. Then as stated, $\mathscr{C}_b^s=B_{\infty,\infty}^s$ whenever $s\in\mathbb{R}^+$, and $\mathscr{C}_b^s=\mathcal{C}_b^s$ whenever $s\in\mathbb{R}^+\setminus\mathbb{Z}$.

We note another useful estimate that is easy to verify by definition (See (1.1) of \cite{debussche2013existence}):
for $f\in C^n(\mathbb{R}^d)$ and all $h\in\mathbb{R}^d$,
\begin{equation}\label{easyinequality2}
\|\Delta_h^n f\|_{L^1}\leq C_n |h|^n\|D^nf\|_{L^1}.
\end{equation}

\subsection{The first estimates}\label{section2.2}

For the first step we improve the regularity of the coefficients.

\begin{theorem}[Theorem 5.1 of \cite{romito2018simple}]
\label{theorem3.1}
Consider the SDE $$dX_t=a(t,X_t)dt+dW_t$$ with initial value $x\in\mathbb{R}^d$ and assume that $a\in L^\infty([0,T]\times\mathbb{R}^d)$. Then the solution has a density $p_x(t)$ for every $t>0$, with the following estimate: for every $\gamma<1$ and $e_\gamma>\gamma$, 
\begin{equation}\label{estimatedensity}\sup_{t\in[0,T]}\sup_{x\in\mathbb{R}^d}(1\wedge t)^{e_\gamma}\|p_x(t)\|_{B_{1,\infty}^\gamma}<\infty.\end{equation}
\end{theorem}

Set the drift $a(t,x):=\langle b(t,x-\cdot),\mu_t\rangle$. Since $b\in L^\infty$, application of the previous theorem shows that for each $t>0$, the law of the SDE has a density in $B_{1,\infty}^\gamma$, that is, $\mu_t$ has a density $p_t\in B_{1,\infty}^\gamma.$ Since $b(t,\cdot)\in L^\infty\subset B_{\infty,\infty}^0(\mathbb{R}^d)$ for each $t>0$, by the convolution inequality \eqref{convolutioninequ} we have that $a(t,x):=\langle b(t,x-\cdot),\mu_t\rangle\in B_{\infty,\infty}^\gamma(\mathbb{R}^d)$ as a function in $x$, for each $t\in[0,T]$ and $\gamma<1$. That is, the drift of the McKean-Vlasov SDE is $\gamma$-Hölder continuous in space for each $\gamma\in(0,1)$, even if the interaction kernel $b$ is only bounded measurable. This sheds some light on how the convolution term $\langle b(t,x-\cdot),\mu_t\rangle$ makes the SDE smoother in law. 

We remark that Theorem \ref{theorem3.1} is true under the following minor modification: for any initial law $\nu_0\in\mathcal{P}(\mathbb{R}^d)$, the SDE has a unique weak solution, whose law $\nu_t$ at each time $t$ is, by Markov property, given by $\nu_t=\int p_x(t)\mu_0(dx).$ This implies $\nu_t$ has a density which we denote by $p_{\mu_0}(t)$, and a direct application of Jensen's inequality implies that $p_{\mu_0}(t)$ satisfies the same estimate as \eqref{estimatedensity}. Therefore we conclude with the following corollary:

\begin{corollary}\label{corollary2.2}
Consider the McKean-Vlasov SDE 
$$dX_t=\langle b(t,X_t-\cdot),\mu_t\rangle dt+dW_t\quad\operatorname{Law}(X_t)=\mu_t$$
with initial law $\mu_0\in\mathcal{P}(\mathbb{R}^d)$. Assume the interaction $b\in L^\infty([0,T]\times\mathbb{R}^d).$ Then the law $\mu_t$ at time $t$ has a density $p_{\mu_0}(t)$ with respect to the Lebesgue measure. Moreover we have the following estimate: for each $\gamma\in(0,1)$,
there exists a decreasing function $\ell_\gamma:(0,T]\to\mathbb{R}^+$, such that 
$$\sup_{t\geq t_0} \|p_{\mu_0}(t)\|_{B_{1,\infty}^\gamma}<\ell(t_0)<\infty,\quad \text{for each }T>t_0>0,$$
$$\sup_{t\geq t_0} \left\|\int  b(t,\cdot-y)p_{\mu_0}(t)(dy)\right\|_{B_{\infty,\infty}^\gamma}<\ell(t_0)<\infty,\quad \text{for each }T>t_0>0.$$
\end{corollary}

\subsection{A step forward}\label{section2.3}

Having proved the drift of the McKean-Vlasov SDE is Hölder when $b\in L^\infty$, we can iterate the above procedure and prove that the law $p_{\mu_0}(t)$ lies in some higher order Besov spaces. The idea of the next few iterations are essentially outlined in Section 2.2 of \cite{romito2018simple}, yet we need to make several key modifications to make the argument work in our case. For this reason, we give the complete details.

The major technical lemma is the following:
\begin{lemma}[Lemma A.1 of \cite{romito2018simple}]\label{lemma4.1} Given a $\mathbb{R}^d$-valued variable $X$, assume for an exponent $\alpha>0$, an index $s>\alpha$ and an integer $m>s$, there exists $K>0$ such that for every $\phi\in \mathscr{C}_b^\alpha(\mathbb{R}^d)$ and $|h|\leq 1$, 
$$\mathbb{E}[|\Delta_h^m\phi(X)]|\leq K|h|^s\|\phi\|_{\mathscr{C}_b^\alpha},$$
then a density $f_X$ exists for $X$, $f_X\in B_{1,\infty}^{s-\alpha}(\mathbb{R}^d)$, and 
$$\|f_X\|_{B_{1,\infty}^{s-\alpha}}\lesssim (1+K).$$
\end{lemma}

Now we prove $p_{\mu_0}(t)$ lies in higher order Besov spaces. We modify Section 2.2 of \cite{romito2018simple} to cover time dependent drifts with $\mathcal{C}_b^\alpha$-norms not uniformly bounded on $[0,T]$. The reason for doing this modification is that the last estimate in Corollary \ref{corollary2.2} is not uniform in time.

Consider the toy-model SDE $$dX_t=a(t,X_t)dt+dW_t,$$ where $a\in L^\infty([0,\infty)\times\mathbb{R}^d)$ and such that for each $\gamma\in(0,1)$ there exists a positive decreasing function  $\ell_\gamma:(0,\infty]\to\mathbb{R}^+$ satisfying
\begin{equation}\label{firstfirstfirst}
\sup_{s\geq t_0}\|a(s,\cdot)\|_{B_{\infty,\infty}^\gamma}\leq \ell_\gamma(t_0)<\infty,\quad \text{ for each } t_0>0.\end{equation}

Now fix some $t\in(0,T]$ and $\delta\in(0,t)$, choose $\epsilon\in (0,(t-\delta)\wedge 1))$ such that \begin{equation}\label{lowrequirement}t-\epsilon>\delta.\end{equation}

To estimate the density of $X_t$ at time $t$, we linearize the system at time $t-\epsilon$, as follows:
define the linearization process
\begin{equation}\label{linearizationprocess}Y_s^\epsilon=X_{t-\epsilon}+\int_{t-\epsilon}^s a(r,X_{t-\epsilon})dr+W_s-W_{t-\epsilon},\quad s\geq t-\epsilon.\end{equation}
Since $$X_s=X_{t-\epsilon}+\int_{t-\epsilon}^s a(r,X_r)dr+W_s-W_{t-\epsilon},$$
a direct computation shows 
\begin{equation}\label{AEestimate}\begin{aligned}
\mathbb{E}[|X_t-Y_t^\epsilon|]&\leq\int_{t-\epsilon}^t |a(r,X_r)-a(r,X_{t-\epsilon})|dr\\&\leq \ell_\gamma(t-\epsilon)\int_{t-\epsilon}^t |X_r-X_{t-\epsilon}|^\gamma ds\\&\leq C\ell_\gamma(\delta)\epsilon^{1+\frac{\gamma}{2}},\end{aligned}\end{equation}
where the second inequality follows from Hölder continuity of $a(r,\cdot)$ and \eqref{firstfirstfirst},  and the last inequality follows from an elementary SDE estimate (see Lemma \ref{elementaryproperty}).

Now conditioning on $X_{t-\epsilon}=y$, $Y_t^\epsilon$ is a shifted Gaussian variable with variance $\epsilon$, so its distribution can be computed explicitly.

By lemma \ref{lemma4.1}, regularity estimates of the law of $X_t$ shall be estimated via evaluating $\mathbb{E}[\Delta_h^m\phi(X_t)]$,  By triangle inequality, 
$$\left|\mathbb{E}[\Delta_h^m\phi(X_t)]\right|=\left|\mathbb{E}[\Delta_h^m\phi(X_t)-\Delta_h^m\phi(Y^\epsilon_t)]\right|+\left|\mathbb{E}[\Delta_h^m\phi(Y^\epsilon_t)]\right|.$$ Following \cite{romito2018simple}, we name the first term on the right hand side  by \textbf{AE}, the approximation error; and name the second term on the right hand side \textbf{PE}, the probability estimate. 

We start with the evaluation of \textbf{PE} via conditioning $\mathbb{E}[\Delta_h^m\phi(Y_t^\epsilon)]=\mathbb{E}[\mathbb{E}[\Delta_h^m\phi(Y_t^\epsilon)]\mathcal{F}_{t-\epsilon}].$ Conditioning on $Y_{t-\epsilon}=y$, $Y_t^\epsilon$ is a Brownian motion with covariance $\epsilon$ and a constant drift, so by \eqref{easyinequality2} and the explicit form of Gaussian density, 
$$|\mathbb{E}[\Delta_h^m\phi(Y_t^\epsilon)]|_{Y_{t-\epsilon}=y}|\leq \|\phi\|_\infty(|h|/\sqrt{\epsilon})^m,\quad y\in\mathbb{R}^d,$$
and taking expectation,  
\begin{equation}\label{PEestimate}\mathbb{E}[\Delta_h^m\phi(Y^\epsilon_t)]\leq \|\phi\|_\infty(|h|/\sqrt{\epsilon})^m.\end{equation}

To evaluate \textbf{AE}$:=\mathbb{E}[\Delta_ h^m\phi(X_t)-\Delta_h^m\phi(Y_t^\epsilon)]$, 
we choose $\phi\in\mathscr{C}_b^\alpha,$ then for some constant $C_m>0$, we have $\|\Delta_h^m\phi\|_{\mathscr{C}_b^\alpha}\leq C_m \|\phi\|_{\mathscr{C}_b^\alpha},$ such that  $$\left|\mathbb{E}[\Delta_ h^m\phi(X_t)-\Delta_h^m\phi(Y_t^\epsilon)]\right|\leq C_m\|\phi\|_{\mathscr{C}_b^\alpha}\mathbb{E}[|X_t-Y_t^\epsilon|^\alpha].$$
Combined with estimates \eqref{AEestimate} and\eqref{PEestimate}, we have 
\begin{equation}\label{modelcase}\left|\mathbb{E}[\Delta_h^m\phi(X_t)]\right|\leq \ell_\gamma(t-\epsilon)\|\phi\|_{\mathscr{C}_b^\alpha}(\epsilon^{\alpha(1+\frac{\gamma}{2})}+(\frac{|h|}{\sqrt{\epsilon}})^m).\end{equation}
Now we choose an appropriate $\epsilon$ to optimize the right hand side. First ignoring the restriction $\epsilon<t-\delta$, then the optimal choice is \begin{equation}\label{choiceepsilon}\epsilon=|h|^\frac{m}{\alpha(1+\frac{\gamma}{2})+\frac{m}{2}},\end{equation} and we have 
$$\left|\mathbb{E}[\Delta_h^m\phi(X_t)]\right|\leq\|\phi\|_{\mathscr{C}_b^\alpha}|h|^\frac{m\alpha(2+\gamma)}{m+\alpha(2+\gamma)} $$
For $m$ large, the exponent in $h$ is $\alpha(2+\gamma),$ by Lemma \ref{lemma4.1}, the law $p_t$ lies in $B_{1,\infty}^{\alpha+\alpha\gamma}.$ Since $\alpha$ can take arbitary value in $(0,1)$, we conclude that $p_t$ is in $B_{1,\infty}^{1+\gamma-}$, that is, the exponent can be anything smaller than $1+\gamma$. 

In the aforementioned estimate we have chosen $\epsilon$ a power of $h$ with $h$ ranging over $(0,1]$, so the estimate is valid whenever $t\geq 1+\delta$ for some prescribed $\delta>0$. In the $t\in(0,1)$ case, a judicious choice of $\epsilon\in(0,1)$ and accordingly a careful choice of $\alpha\in(0,1)$ will be necessary. We summarize the choices in the following proposition, which is a minor generalization of Proposition 2.2 of \cite{romito2018simple}. The estimate in \eqref{modelcase} corresponds to setting $a_0=1+\gamma$. We stress that we require $\epsilon\leq\frac{t}{2}$ in the Proposition, which corresponds to taking $\delta=\frac{t}{2}$ in \eqref{lowrequirement}.

\begin{proposition}[see \cite{romito2018simple},  Proposition 2.2]\label{technicalproposition} Consider a solution $X_t$ to the SDE
$$dX_t=a(t,X_t)dt+dW_t$$ with initial distribution $\mu_0\in\mathcal{P}(\mathbb{R}^d).$

Suppose there are numbers $a_0>0$ and a \textbf{decreasing} function $\ell:t\in[0,\infty)\to\mathbb{R}^+$ such that 
$$|\mathbb{E}[\Delta_h^m\phi(X_t)]|\leq \ell(t)\left(e^{\frac{\alpha}{2}(1+a_0)}+(\frac{|h|}{\sqrt{\epsilon}})^m\right)\|\phi\|_{\mathscr{C}_b^\alpha}$$
holds for all $\epsilon\leq 1$, $\epsilon\leq \frac{t}{2}$, all $\alpha\in(0,1)$ and each $\phi\in \mathscr{C}_b^\alpha$, then for any $a\in(0,a_0)$, there exists a decreasing function $\ell_a:(0,T]\to \mathbb{R}^+$ such that $P_x(t)\in B_{1,\infty}^a$ for any $t>0$ and 
\begin{equation}
\sup_{s\geq t_0}\|p_x(s)\|_{B_{1,\infty}^a}\leq \ell_a(t_0)<\infty,\quad\text{ for each } T> t_0>0.
\end{equation}

\end{proposition}

\begin{remark}
Suppose that the function $\ell$ in the assumption is constant, then the function $\ell_a(t)$ in the conclusion is of the form $(1\wedge t)^{-c}$ for some constant $c=c(a,a_0)>0$, see \cite{romito2018simple}. In our induction procedure, $\ell$ is in general unbounded, and the precise asymptotic of $\ell_a$ is not needed in the proof, so we choose to state the proposition in this \textit{qualitative} form. 
\end{remark}

The proof of this proposition is similar to that in \cite{romito2018simple}. We give a sketch of proof in Appendix \ref{appendixA} for sake of completeness. The idea goes as follows: when $t$ is large, $\epsilon$ is chosen to be a suitable power of $h$, as did in \eqref{choiceepsilon}; when $t$ is small, we simply choose $\epsilon=\frac{t}{2}.$ It remains to optimize the various exponents and then conclude via applying Lemma \ref{lemma4.1}.

\begin{corollary}\label{corostep1}
In the setting of Corollary \ref{corollary2.2}, for each $\gamma\in(0,1)$ and each $t\in[0,T]$, the density $\mu_t$ has a density $p_t$ that lies in $B_{1,\infty}^{1+\gamma}(\mathbb{R}^d)$ and the drift $x\to \langle b(t,x-\cdot),\mu_t\rangle$ is $1+\gamma$-Hölder continuous. More precisely, for each $\gamma\in(0,1)$ we can find a decreasing function $\ell_{1+\gamma}:[0,T]\to\mathbb{R}^+$ such that for each $t_0\in(0,T],$
$$\sup_{s\geq t_0>0}\|p_s\|_{B_{1,\infty}^{1+\gamma}}\leq \ell_{1+\gamma}(t_0)<\infty,$$
and
$$\sup_{s\geq t_0>0} \|x\mapsto \langle b(s,x-\cdot),\mu_s\rangle\|_{B_{\infty,\infty}^{1+\gamma}}\leq \ell_{1+\gamma}(t_0)<\infty.$$
\end{corollary}

Let us summarize what we have proved. The starting point is the modest hypothesis $b\in L^\infty$, then we show that the drift is $\gamma$-Hölder continuous for any $\gamma\in(0,1)$, from which we show the drift is $\gamma$-Hölder continuous for any $\gamma\in(1,2)$. The same takes place when we consider the marginal law $\mu_t$: we first show its density $p_t$ lies in $B_{1,\infty}^\gamma$ for $\gamma\in(0,1)$, then show it lies in $B_{1,\infty}^\gamma$ for $\gamma\in(1,2)$. We will make one step further in the next subsection.

\subsection{Second order regularity, limitations of the present approach}\label{section2.4}

We are back to our model case $dX_t=a(t,X_t)dt+dW_t$ with $a\in L^\infty$ and $a(t,\cdot)\in \mathcal{C}_b^{1+\gamma}$ for all $\gamma\in(0,1).$ Assume that for each $\gamma\in(0,1)$ there exists a positive decreasing function  $\ell_{1+\gamma}:(0,\infty]\to\mathbb{R}^+$ satisfying
\begin{equation}\label{secondsecondsecond}
\sup_{s\geq t_0}\|a(s,\cdot)\|_{B_{\infty,\infty}^{1+\gamma}}\leq \ell_{1+\gamma}(t_0)<\infty,\quad \text{ for each } t_0>0.\end{equation} We introduce a higher order linearization 
$$Y_s^\epsilon=X_{t-\epsilon}+\int_{t-\epsilon}^s A_rdr+W_s-W_{t-\epsilon},\quad s\geq t-\epsilon$$
with $A_r=a(r,X_{t-\epsilon})+Da(r,X_{t-\epsilon})(W_r-W_{t-\epsilon}).$

We evaluate $a(s,X_s)-A_s$ as in Section 2.2 of \cite{romito2018simple}: for any $s>t$, 
\begin{equation}\begin{aligned}
    a(s,X_s)-A_s=&a(s,X_s)-a(s,X_{t-\epsilon})-Da(s,X_{t-\epsilon})(W_s-W_{t-\epsilon})\\
    =&Da(s,X_{t-\epsilon})(X_s-X_{t-\epsilon})+O(|X_s-X_{t-\epsilon}|^{1+\gamma})-Da(s,X_{t-\epsilon})(W_s-W_{t-\epsilon})\\=&Da(s,X_{t-\epsilon})\int_{t-\epsilon}^s a(r,X_r)dr+O(|X_s-X_{t-\epsilon}|^{1+\gamma})\\\leq&\ell_{1+\gamma}(t-\epsilon)(s-t+\epsilon)^{\frac{1+\gamma}{2}}.    \end{aligned}
\end{equation}

In the last line we have used the following elementary property:
\begin{lemma}\label{elementaryproperty}
Consider the SDE $dX_t=a(t,X_t)dt+dW_t$ with $a\in L^\infty,$ then for each $m\geq 1$, there exists a constant $C_m>0$ such that 
$$\mathbb{E}[|X_s-X_t|^m]\leq C_m |s-t|^{\frac{m}{2}},\quad s\geq 0,t\geq 0, |s-t|<1.$$
\end{lemma}

\begin{proof}
For $t>s$, $X_t-X_s=\int_s^t a(r,X_r)dr+W_t-W_s.$  Note that $W_t-W_s$ has the same law as $W_{t-s}$. We use explicit form of the Gaussian density, the elementary inequality $(u+v)^m\leq C_m(u^m+v^m)$, $a\in L^\infty$, and $|s-t|<1$ to conclude.
\end{proof}

Now we compute the approximation error, \textbf{AE}. For any $\alpha\in(0,1)$,
$$\mathbb{E}[|X_t-Y_t^\epsilon|^\alpha]\leq \mathbb{E}[|X_t-Y_t^\epsilon|]^\alpha\leq\ell_{1+\gamma} \left(\int_{t-\epsilon}^t (s-t+\epsilon)^{\frac{1+\gamma}{2}}ds\right)^\alpha=\ell_{1+\gamma}\epsilon^{\frac{\alpha}{2}(3+\gamma)}.$$
Therefore for any $\phi\in\mathscr{C}_b^\alpha$ and $m$ sufficiently large, the approximation error \textbf{AE} is
$$\mathbb{E}[|\Delta_h^m\phi(X_t)-\Delta_h^m\phi(Y_t^\epsilon)|]\leq\ell_{1+\gamma}(t-\epsilon) \|\phi\|_{\mathscr{C}_b^\alpha}\epsilon^{\frac{\alpha}{2}(3+\gamma)}.$$

For the probabilistic estimate \textbf{PE}, conditioning on $X_{t-\epsilon}=y$, we have
$$Y_t^\epsilon=X_{t-\epsilon}+\int_{t-\epsilon}^t a(r,y)dr+\int_{t-\epsilon}^t Da(r,y)\widetilde{W}_{r-t+\epsilon}dr+\widetilde{W}_{\epsilon},$$
where  $\widetilde{W}_r:=W_{t-\epsilon+r}-W_{t-\epsilon}$ has the law of a standard Brownian motion for $r\geq 0.$

Conditioning on $Y_{t-\epsilon}=y$, the law of $Y_t^\epsilon$ is Gaussian with variance
$$ \operatorname{Var}(Y_t^\epsilon|_{X_{t-\epsilon}=y})= \int_0^\epsilon \left[I_d+\int_s^\epsilon Da(r+t-\epsilon,y)dr\right]\left[I_d+\int_s^\epsilon Da(r+t-\epsilon,y)dr\right]^Tds.$$
By the assumption \eqref{secondsecondsecond}, $Da$ is locally bounded. Find two constants $0<c_1<<1<<c_2$, then for each $t>0$ we may find some $h_t\in(0,1)$ such that   $$c_1\epsilon I_d\leq \operatorname{Var}(Y_t^\epsilon|_{X_{t-\epsilon}=y})\leq c_2\epsilon I_d,\quad\text{ given } 0<\epsilon\leq h_tt<t.$$ 

Therefore we have the same probability estimate  \textbf{PE} as in the previous steps:
$$\mathbb{E}[|\Delta_h^m\phi(Y_t^\epsilon)|]\leq\|\phi\|_{\mathscr{C}_b^\alpha} (\frac{h}{\sqrt{\epsilon}})^m,\quad \text{ given } 0<\epsilon\leq h_tt<t.$$

Combining the above estimates, we have (compare with \eqref{modelcase})
\begin{equation}\label{nonmodelcase}\left|\mathbb{E}[\Delta_h^m\phi(X_t)]\right|\leq \ell_{1+\gamma}(t-\epsilon)\|\phi\|_{\mathscr{C}_b^\alpha}(\epsilon^{\frac{\alpha}{2}(3+\gamma)}+(\frac{|h|}{\sqrt{\epsilon}})^m),\quad 0<\epsilon\leq h_t t<t.\end{equation}

Due to the additional restriction $\epsilon\leq h_t t$, some modifications will be made when applying Proposition \ref{technicalproposition}. The details of these modifications are contained in Appendix \ref{appendixA}. Nonetheless, we can still upgrade the regularity in Corollary \ref{corostep1} and conclude with:

\begin{corollary}\label{finalcorollary}
In the setting of Corollary \ref{corollary2.2}, for each $\gamma\in(0,1)$ and each $t\in[0,T]$, the density $\mu_t$ has a density $p_t$ that lies in $B_{1,\infty}^{2+\gamma}$ and the drift $x\to \langle b(t,x-\cdot),\mu_t\rangle$ is $2+\gamma$-Hölder continuous. More precisely, for each $\gamma\in(0,1)$ we can find a decreasing function $\ell_{2+\gamma}:[0,T]\to\mathbb{R}^+$ such that for each $t_0\in(0,T],$
$$\sup_{s\geq t_0>0}\|p_s\|_{B_{1,\infty}^{2+\gamma}}\leq \ell_{2+\gamma}(t_0)<\infty$$
and
$$\sup_{s\geq t_0>0} \|x\mapsto \langle b(s,x-\cdot),\mu_s\rangle\|_{B_{\infty,\infty}^{2+\gamma}}\leq \ell_{2+\gamma}(t_0)<\infty.$$
\end{corollary}

We have now proved the first part of Theorem \ref{theorem1}. In dimension one, the Sobolev embedding $B_{1,\infty}^s\hookrightarrow B_{\infty,\infty}^{s-1}$ for $s\in\mathbb{R}$ implies that the density $p_t$ is $1+\gamma$- Hölder continuous for any $\gamma\in(0,1)$. This proves Theorem \ref{theorem1} completely in the one dimensional case, and shows that $p_t$ already more regular than the density of the SDE introduced in Appendix \ref{appendixA}, which is not of McKean-Vlasov type.

The remaining parts of Theorem \ref{theorem1} will be proved in Section \ref{section4}.

\begin{remark} \label{remark1}
One may hope to iterate the present argument a number of times, utilizing stochastic Taylor expansions $Y_t^\epsilon$ of higher orders, and eventually prove that the density $p_t$ is smooth. However a technical problem arises: when we do higher order expansions, we encounter iterated (and much higher orders of) stochastic integrals of the form $\int_0^t\int_0^s W_rdW_rdW_s,$ which is no longer Gaussian. Consequently the probability density function of $Y_t^\epsilon$ cannot be easily analyzed, and we are not sure how to estimate the term \textbf{PE}. 

One alternative is to use the standard Malliavin calculus toolbox, as we pursue in the next section. The price to pay is we will impose additional regularity assumptions on $b$.
\end{remark}

\subsection{A route towards smoothness}\label{section2.5}

Utilizing stochastic Taylor expansions of higher order (see for example \cite{kloeden1992stochastic}, Section 10), we can expand $X_t$ at $t-\epsilon$ to much higher orders. Then $Y_t$, the approximating process, involves multiple stochastic integrals lying in several Wiener chaos.

To control the probability estimate \textbf{PE} without knowing the probability density function, one may use Malliavin calculus (the notions and techniques of Malliavin calculus are outlined in Section \ref{section3}), noting that $Y_t$ is Malliavin smooth. Then the key computation is contained in Lemma \ref{lemma3.6}, which involves estimating the $L^p$-norms of the Malliavin matrix $\mathscr{D}Y_t$ of $Y_t$, as well as estimating $p$-th norms of the inverse of the Malliavin matrix. We should show the former has order $O(1)$, and the latter has order $O(\epsilon)$. The most difficult part is the norms of the inverse Malliavin matrix. As we can only estimate the $p$-th norms of the Malliavin matrix, we can only estimate the norm of the inverse following the steps given in Appendix \ref{momentappendix}, which does not seem to give sharp results.

\section{High regularity regime: Malliavin Calculus}\label{section3}
In this section we prove Theorem \ref{theorem2}.
We use Malliavin calculus to upgrade regularity of the density $p_t$ via a self-improvement procedure: assuming $b$ has some regularity, then $\langle b(x-\cdot,\mu_t\rangle$ has a higher regularity than $p_t$ itself. Now we apply Malliavin calculus given this better regularity, and this improves regularity of $p_t$ a further step.

In this procedure, a sharp control of the number of times the process is Malliavin differentiable is critical to our analysis. The sharpest result we are aware of is \cite{banos2016malliavin}.

Section \ref{section3.4} is the central part of our analysis, whose results are then applied in Section \ref{section3.5} to McKean-Vlasov SDEs. Section  \ref{section3.1} establishes the Malliavin calculus framework, then in section \ref{section3.2} and \ref{section3.3} essential technical preparations are set up.

\subsection{The Malliavin calculus framework}
\label{section3.1}
We recall briefly the framework of Malliavin calculus in this section, more details can be found in \cite{nualart2006malliavin} or \cite{hairer2021introduction}. 

We are given a real separable Hilbert space $H$ (in this paper we take $H=L^2(0,T];\mathbb{R}^d)$) and a Gaussian white noise process $W=\{W(h),h\in H\}$ defined on a complete probability space $(\Omega,\mathcal{F},P)$ with the filtration generated by $W$.

Let $\mathcal{S}\subset L^2(\Omega,P)$ be the subspace of random variables $X$ that have the simple form 
$$X=f(W(h_1),\cdots,W(h_N)),$$  for some $N\in\mathbb{N}_+$, some $\mathcal{C}_0^\infty$ function $f:\mathbb{R}^N\to\mathbb{R}$, and some elements $h_i\in H$. On $\mathcal{S}$ we define a differential operator $\mathscr{D}$ which sends $X$ to a $H$-valued random variable $\mathscr{D}X$:
$$\mathscr{D}X=\sum_{k=1}^N \partial_k f(W(h_1),\cdots,W(h_N)) h_k.$$
We then inductively define $\mathscr{D}^k$ for each $k\geq 1$, which sends a random variable $X\in\mathcal{S}$ to an $H^{\otimes k}$-valued random variable $\mathscr{D}^kX$. 

For each $k\geq 1$ and $p\geq 1$, $\mathscr{D}^k$ is a closable operator from $L^p(\Omega)$ to $L^p(\Omega;H^{\otimes k})$. We denote its closed extension by $\mathscr{D}^k$ as well, whose domain contains the space $\mathbb{D}^{k,p}$ defined as follows: for each integer $k\geq 1$ and real number $p\geq 1$, $\mathbb{D}^{k,p}$ is the closure of $\mathcal{S}$ with respect to the norm $$\|X\|_{\mathbb{D}^{k,p}}:=\|X\|_{L^p(\Omega)}+\sum_{i=1}^k\|\mathscr{D}^{i}X\|_{L^p(\Omega;H^{\otimes i})}.$$
More generally, for any real separable Hilbert space $V$, consider the space of $V$-valued simple functions
$$\mathcal{S}_V:=\{u=\sum_{j=1}^N X_j h_j:X_j\in \mathcal{S},h_j\in V,N\geq 1\},$$
and define analogously (using the same symbol $\mathscr{D}$ for differentiation):
$$
\mathscr{D}^ku=\sum_{j=1}^n\mathscr{D}^kX_j\otimes h_j,\quad k\in\mathbb{N}_+.
$$ The operator $\mathscr{D}$ (and its powers) is closable as well. Then for each $k$ and $p$ we define analogously $\mathbb{D}^{k,p}(H)$ as the closure of $\mathcal{S}_V$ with respect to the norm
$$\|X\|_{\mathbb{D}^{k,p}(V)}:=\|X\|_{_{L^p(\Omega;V)}}+\sum_{i=1}^k\|\mathscr{D}^{i}X\|_{L^p(\Omega;H^{\otimes i}\otimes V)}.$$

The divergence operator $\delta$ is defined as the dual of $\mathscr{D}$. That is, $\delta$ is an unbounded operator from $L^2(\Omega;H)$ to $L^2(\Omega).$ Let $\operatorname{Dom}(\delta)$ denote the domain of $\delta$, it is not hard to prove (see for example Theorem 3.6 of \cite{hairer2021introduction}) that $\mathbb{D}^{1,2}(H)\subset \operatorname{Dom}(\delta)$ and $\delta$ is continuous from $\mathbb{D}^{1,2}(H)$ to $L^2(\Omega)$: for $u\in \mathbb{D}^{1,2}(H),$ $\|\delta u\|_{L^2}\leq \|u\|_{\mathbb{D}^{1,2}(H)}$. Indeed, this inclusion holds for all indices $k$ and $p$, as given by the following celebrated Meyer's inequality:

\begin{proposition}[\cite{nualart2006malliavin},Proposition 1.5.7] \label{meyerinequality}
for every $k\geq 1$ and $p>1$, the operator $\delta$ maps $\mathbb{D}^{k,p}(H)$ continuously into $\mathbb{D}^{k-1,p}$. That is, for $u\in \mathbb{D}^{k,p}(H),$
$$\|\delta(u)\|_{\mathbb{D}^{k-1,p}}\leq C_{k,p}\|u\|_{\mathbb{D}^{k,p}(H)}.$$ In particular, $\mathbb{D}^{k,p}(H)\subset \operatorname{Dom}(\delta).$

\end{proposition}

\subsubsection{The context of function spaces}

The special case $H=L^2([0,T];\mathbb{R}^d)$ is of particular interest. We use the identification $L^2(\Omega;H):=L^2(\Omega\times [0,T])$. Then the Malliavin derivative of the process
$X=f\left(\int_0^T h_1(u)dW_u,\cdots,\int_0^T h_n(u)dW_u\right)$
can be regarded as a process $\{\mathscr{D}_tF\}_{t\in[0,T]}$ in $L^2(\Omega\times [0,T];\mathbb{R}^d)$ defined as
$$\mathscr{D}_tX=\sum_{i=1}^n \frac{\partial}{\partial x_i}f\left(\int_0^T h_1(u)dW_u,\cdots,\int_0^T h_n(u)dW_u\right)h_i(t).$$

With this identification, the $\mathbb{D}^{k,p}$ norm of $X$ can be reformulated as
$$\|X\|_{\mathbb{D}^{k,p}}=\mathbb{E}[|F|^p]^{1/p}+\sum_{i=1}^k \mathbb{E}\left[\int_0^T\cdots \int_0^T\|\mathscr{D}_{t_1}\cdots \mathscr{D}_{t_i}F\|^p dt_1\cdots dt_i\right]^{1/p}.$$

Malliavin derivative satisfies the chain rule: let $\varphi:\mathbb{R}^d\to \mathbb{R}^d$ be continuously differentiable. Given a random variable $X=(X_1,\cdots,X_d)$ with each component $X_i\in\mathbb{D}^{1,2}$, we have $\varphi(X)\in\mathbb{D}^{1,2}$ with 
$$\mathscr{D}\varphi(X)=\sum_i \partial_i\varphi(X) \mathscr{D}X_i.$$

\subsubsection{Malliavin differential of a diffusion process}
We specialize to the setting of diffusion processes. Let $a:[0,T]\times\mathbb{R}^d\to\mathbb{R}^d$ be smooth and of linear growth, consider the SDE $X_t=x+\int_0^t b(u,X_u)du+W_t$. Fix $s_1,t\in[0,T]$ and $s_1<t$, by the chain rule of Malliavin calculus we obtain
$$\mathscr{D}_{s_1}X_t=\mathcal{I}_d+\int_{s_1}^t a'(u,X_u)\mathscr{D}_{s_1}X_u du.$$
As suggested in the proof of \cite{banos2016malliavin} Theorem 3.4, we may use Picard iteration to obtain the following series expansion
\begin{equation}\label{theseriesexpansion}
    \mathscr{D}_{s_1}X_t=\mathcal{I}_d+\sum_{m\geq 1}\int_{\Lambda_m(s_1,t)}a'(u_1,X_{u_1})\cdots a'(u_m,X_{u_m})du_1\cdots du_m,
\end{equation}
where for each $m\in\mathbb{N}_+$ define  $\Lambda_m(s,t):=\{(u_1,\cdots,u_m)\in [0,T]^m:s<u_1<\cdots<u_m<t\}.$

\subsection{Malliavin differentiability with optimal index}
\label{section3.2}

We will use the following result in \cite{banos2016malliavin}, which is not covered by classical monographs on Malliavin calculus.

\begin{proposition}[Theorem 3.4 and Theorem 3.7 of \cite{banos2016malliavin}]\label{propositionuseful}
Consider the SDE
$$dX_t=a(t,X_t)dt+dW_t,\quad X_0=x$$
with $a(t,\cdot)$ differentiable in $x$ for each $t$, such that  $a(t,x),Da(t,x)\in L^\infty([0,T]\times\mathbb{R}^d).$  Then 
\begin{enumerate}
    \item   $X_t\in\mathbb{D}^{2,p}(\Omega)$ for each $p\geq 1$. 

\item In general, assume that for some $k\in\mathbb{N}_+$, $$a(t,x),Da(t,x),D^2a(t,x),\cdots,D^ka(t,x)\in L^\infty([0,T]\times \mathbb{R}^d).$$ Then $X_t\in\mathbb{D}^{k+1,p}(\Omega)$ for each $p\geq 1$.
\end{enumerate}
\end{proposition}

As the result of this Proposition is crucial for our application, we give a sketch of proof to show why this proposition is true. Clearly, part (2) follows from induction once we have proved part (1). To prove part (1), we find a sequence of smooth functions $a^n$ approximating $a$, and consider $X^n$ solution to the SDE 
$$dX^n_t=a^n(t,X^n_t)dt+dW_t.$$ 

To prove $X_t\in\mathbb{D}^{2,p}(\Omega),$ it suffices to prove $X_t^n\to X_t$ in $L^p(\Omega)$ and $$\sup_{n\in\mathbb{N}_+}\|X_t^n\|_{\mathbb{D}^{2,p}(\Omega)}<\infty.$$ 

The first claim is easy to verify and we omit it. To prove the second, we differentiate $X_t^n$ twice. Then in the expression of $D_{s_2}D_{s_1}X_t^n$, we have terms involving both $Da^n$ and $D^2a^n$. Notice that we have assumed $Da\in L^\infty$ but not even assuming $D^2a$ exists. Therefore we have to get rid of $D^2a^n$ terms in our estimate. This is facilitated by the following result:
\begin{lemma}[Proposition 3.7 of \cite{menoukeu2013variational} or Proposition 3.3 of \cite{banos2016malliavin}] For a Brownian motion $W$ starting from $z_0\in\mathbb{R}^d$, $b_1,\cdots,b_m$ compactly supported $C^1$ functions, $\alpha_1,\cdots,\alpha_m\in\{0,1\}^d$ multiindex with $|a_i|\leq 1$ for each $i$, we can find a universal constant $C$ such that
$$\left|\mathbb{E}\left[\int_{t_0<t_1<\cdots<t_m<t}(\prod_{i=1}^m D^{\alpha_i}b_i(t_i,B_{t_i}))dt_1\cdots dt_m\right]\right|\leq \frac{C^m\prod_{i=1}^m \|b_i\|_\infty (t-t_0)^{m/2}}{\Gamma(\frac{m}{2}+1)}.$$
\end{lemma}

From this lemma, one immediately sees that the SDE $dX_t=a(t,X_t)dt+dW_t$ with $a\in L^\infty$ is Malliavin differentiable, as proved in \cite{menoukeu2013variational}. Proposition \ref{propositionuseful} generalizes this result to higher order derivatives. It is also noted in \cite{banos2016malliavin}, Example 4.6 that the degree of differentiability in Proposition \ref{propositionuseful} is sharp: for some drift $a$ with $Da\in L^\infty$ and $a$ of linear growth, the corresponding SDE $X_t$ satisfies $X_t\in\mathbb{D}^{2,p}(\Omega)$ but $X_t\notin\mathbb{D}^{3,p}(\Omega)$.

\begin{claim}\label{claims0}
Under the same assumption as Proposition \ref{propositionuseful} (2), for any finite $T>0$ we have the following bound 
$$\sup_{t\in[0,T]}\|X_t\|_{\mathbb{D}^{k+1,p}(\Omega)}\leq C_{p,T}<\infty,\quad p\geq 1,$$
where the constant $C_{p,T}$ depends only on the $L^\infty$ norm of $a,Da,\cdots,D^k a$, $T$ and $p$, and is independent of the initial value $x\in\mathbb{R}^d.$
\end{claim}

To verify this claim one just needs to carefully go through the proof of Proposition \ref{propositionuseful}.

\subsection{The Malliavin matrix is invertible}\label{section3.3}

In this section we study the Malliavin matrix $$\gamma_{X_t}^{ij}:=\langle D_\cdot X_t^{(i)},D_\cdot X_t^{(j)}\rangle_H,\quad i,j=1,\cdots,d.$$

The celebrated Hörmander's theorem (see for example \cite{hairer2021introduction}) says that assuming a parabolic Hörmander's condition on the vector fields of a diffusion process (the iterated Lie brackets of the driving vector fields are non-degenerate at each point), the Malliavin matrix $\gamma_{X_t}^{ij}$ of $X$ is invertible, and has inverse moments of all orders. The invertibility of $\gamma_{X_t}^{ij}$ then implies the diffusion process has a smooth density.

In our (simplified) setting, the nondegenerate additive Brownian noise guarantees invertibility of $\gamma_{X_t}^{ij}$, and its inverse moments of all orders.

The first part of the following proposition is given in Proposition 4.4 of \cite{banos2016malliavin}.

\begin{proposition}\label{claims1}  Let $X_t$ be the SDE in Proposition \ref{propositionuseful}, with $a$ satisfying the same assumption, that is, $a\in L^\infty$, $Da\in L^\infty$ uniformly in $t\in[0,T]$. Then the Malliavin matrix $\gamma_{X_t}$ is almost surely invertible, and 
$$(\det\gamma_{X_t})^{-1}\in \cap_{p\geq 1}L^p(\Omega).$$
Moreover, for each $t_0>0$, we have a uniform in time estimate 
$$\sup_{t\geq t_0}\|(\det\gamma_{X_t})^{-1}\|_{L^p}<C(p,t_0)<\infty,\quad p\geq 1$$
 for some constant $C(p,t_0)$ depending on $p$ and $t_0$. This estimate does not depend on the initial value $x$ of the SDE $X_t$.
\end{proposition}

\begin{proof}
This is essentially the same as Proposition 4.4 of \cite{banos2016malliavin}, the difference is we need to make explicit each quantitative estimate that appears in the proof, as we aim to show that the final estimates are uniform over $t>t_0>0$. A complete proof of this proposition is given in Appendix \ref{momentappendix}.
\end{proof}

 \begin{claim}\label{remark3.4}
 For every $k\in\mathbb{N}_+$ define $\mathcal{S}^k:=\cap_{p\geq 1}\mathbb{D}^{k,p}(\Omega)$. Then $\mathcal{S}^k$ is an algebra by the Leibniz rule and Hölder's inequality. Similarly, for each Hilbert space $K$ define $\mathcal{S}^k(K):=\cap_{p\geq 1}\mathbb{D}^{k,p}(\Omega;K).$ $\mathcal{S}^k(K)$ is also an algebra
 \end{claim}
 
  The following result is also useful.

 \begin{lemma} Assume $X_t$ solves the SDE $dX_t=a(t,X_t)dt+dW_t$ where $a,Da,\cdots,D^ka\in L^\infty,$ i.e., $a$ satisfies the assumption (2) of Proposition \ref{propositionuseful} with index $k$.
Then the inverse of the Malliavin matrix satisfies $\gamma_{X_t}^{-1}\in \mathcal{S}^{k+1}$.
 \end{lemma}
 
 \begin{proof}
The chain rule of $\mathscr{D}$ implies that $\mathscr{D}^{k+1}\gamma_{X_t}^{-1}$ can be written as a polynomial of $\gamma_{X_t}^{-1}$ and $\mathscr{D}^\ell X_t$ for $\ell\leq k+1$. Also by Proposition \ref{propositionuseful} (2), $X_t\in \mathbb{D}^{k+1,p}$ for each $p\geq 1$. The result then follows from Hölder's inequality.
 \end{proof}

The next lemma reformulates Lemma 4.4 of \cite{hairer2021introduction} in the setting where $X$ is not assumed Malliavin smooth. This is essentially the integration by parts formula in Malliavin calculus, see for example Proposition 22 of \cite{bally2011riesz}.

\begin{lemma}\label{lemma3.6}
Let $X$ be an $\mathbb{R}^d$-valued random variable satisfying: for some $k\in\mathbb{N}_+$, $X\in\mathcal{S}^k$, and the Malliavin matrix $\gamma_X$ of $X$ is almost surely invertible with $\gamma_X^{-1}\in \mathcal{S}^k.$ Then for each $Z\in\mathcal{S}^{k-1}$ we can find a $\bar{Z}\in\mathcal{S}^{k-2}$ such that
\begin{equation}
    \label{integralbyparts}
\mathbb{E}[Z\partial_i G(X)]=\mathbb{E}[G(X)\bar{Z}]\end{equation}
holds for each $G\in\mathcal{C}_0^\infty.$
\end{lemma}

\begin{proof}
 By chain rule, $\mathscr{D}G(X)=\sum_j\partial_j G(X)\mathscr{D}X_j.$ Then if we take 
 $$Y_i=\sum_{j=1}^n(\mathscr{D}X_j)\gamma_{X_{ji}}^{-1},$$ we have the identity $\partial_i G(X)=\langle\mathscr{D}G(X),Y_i\rangle.$
 Therefore if we choose $\bar{Z}=\delta(ZY_i),$ the relation
 \eqref{integralbyparts}
 holds. It remains to check $\bar{Z}\in\mathcal{S}^{k-2}$.
 By Claim \ref{remark3.4} and assumptions on $X$, we see that $Y_i\in\mathcal{S}^{k-1}(H)$. Then the assumption $Z\in\mathcal{S}^{k-1}$ implies $ZY_i\in\mathcal{S}^{k-1}(H)$, then Meyer inequality \ref{meyerinequality} implies that $\bar{Z}\in\mathcal{S}^{k-2}$.
\end{proof}

\subsection{Regularity estimates for the density}\label{section3.4}

The following lemma transfers Malliavin differentiability into smoothness of the density. It is inspired by Lemma 4.1 of \cite{hairer2021introduction}, but we only assume differentiability up to order $k$ and state the estimates in terms of Besov norms. This lemma is also similar in spirit to Lemma \ref{lemma4.1}. For duality of Besov spaces, one may refer to Appendix A.1 of \cite{romito2018simple}.

\begin{lemma}\label{lemma3.80}
Let $X$ be an $\mathbb{R}^d$-valued random variable where there exists $C_m$ such that $|\mathbb{E}D^{(m)}G(X)|\leq C_m\|G\|_\infty$ for every $G\in\mathcal{C}_0^\infty$ and $m=1,\cdots,k$. Then the law of $X$ has a density with respect to the Lebesgue measure that lies in the Besov space $B_{1,1}^{k}(\mathbb{R}^d)\hookrightarrow B_{1,\infty}^{k}(\mathbb{R}^d).$
\end{lemma}

\begin{proof} 
Let $\mu$ denote the law of $X$, then the assumption translates to 
$$\left|\int_{\mathbb{R}^d}D^{(m)}G(x)\mu(dx)\right|\leq C_m\|G\|_\infty,\quad m=1,\cdots,k.$$
By the continuous embedding $L^\infty\hookrightarrow B_{\infty,\infty}^0$, duality of $B_{\infty,\infty}^0$ with $B_{1,1}^0$, and denseness of $\mathcal{C}_0^\infty$ in $B_{\infty,\infty}^0$, the assumption implies the distributional derivatives of $\mu$ up to order $k$ lie in $B_{1,1}^0$. Consequently $\mu$ has a density that lies in $B_{1,1}^k$, which embeds continuously in $B_{1,\infty}^k$.
\end{proof}

\begin{theorem} \label{theorem3.9}
Let $X$ be an $\mathbb{R}^d$-valued random variable satisfying: for some integer $k\geq 2$,  $X\in\mathcal{S}^k$ and the inverse of its Malliavin matrix $\gamma_X^{-1}\in\mathcal{S}^k$. Then $X$ has a density $p_X$ with respect to the Lebesgue measure with $p_X\in B_{1,\infty}^{k-1}(\mathbb{R}^d).$ 
\end{theorem}

\begin{proof}
We prove the following claim for each $m=1,2,\cdots,k-1$: for every $Y\in\mathcal{S}^{k-m}$ there exists $Z\in\mathcal{S}^{k-m-1}$ such that 
\begin{equation}
\label{proofbyinduction}    
\mathbb{E}[YD^{(m)}G(X)]=\mathbb{E}[G(X)Z]\end{equation}
for every $G\in\mathcal{C}_0^\infty.$ Then the claim follows by taking $Y=1$, $m=k-1$ and checking that $|\mathbb{E}[G(X)\bar{Z}]|\leq \|G\|_\infty \mathbb{E}[|\bar{Z}|]<\infty. $
The claim \eqref{proofbyinduction} is trivial for $m=0$, and is true for $m=1$ by Lemma \ref{lemma3.6}. Suppose it is true for some $\ell<k-1$, then 
\begin{equation}\label{proofbyinduction12}
\mathbb{E}[YD^{(\ell+1)}G(X)]=\mathbb{E}[DG(X)Z]=\mathbb{E}[G(X)\bar{Z}]\end{equation}
for some $\bar{Z}\in\mathcal{S}^{k-\ell-2}$ by Lemma \ref{lemma3.6}. Thus \eqref{proofbyinduction} is also true for $m=\ell+1$, and the proof follows by induction.
\end{proof}

\begin{claim}\label{claims2}
With the same assumption as in Theorem \ref{theorem3.9}, the Besov norm  $\|p_X\|_{B_{1,\infty}^{k-1}}$ can be bound by a constant $C(X,\gamma_X)$ depending only on $(\|X\|_{\mathbb{D}^{k,p}})_{p\geq 1}$ and $(\|(\det\gamma_X)^{-1}\|_{L^p})_{p\geq 1}.$
\end{claim}

\begin{proof}
By duality and Lemma \ref{lemma3.6}, it suffices to bound $\|\bar{Z}\|_{L^1}$ where $\bar{Z}$ is the random variable appearing in equation \eqref{proofbyinduction12} (with the choice $Y=1$ and $\ell=k-2$). From the proof of Lemma \ref{lemma3.6}, $\bar{Z}$ can be represented as a $k-1$-fold implementation of the divergence operator $\delta$ and in each layer there are occurrences of $\gamma_X^{-1}$ and $\mathscr{D}X$. Thus to bound $\|\bar{Z}\|_{L^1}$ one only needs to use Hoölder's inequality and Meyer's inequality (Proposition \ref{meyerinequality}) $k-1$ times.
\end{proof}

Gathering the various claims in this section, we have the following corollary

\begin{corollary}\label{corollarystep1}
Consider the SDE $dX_t=a(t,X_t)dt+dW_t$, $X_0=x\in\mathbb{R}^d$ with the drift $a$ satisfying $a,Da,\cdots,D^ka\in L^\infty([0,T]\times\mathbb{R}^d)$. Then for each $t>0$, the law of $X_t$ has a density  $p_t\in B_{1,1}^{k}(\mathbb{R}^d)$ with respect to the Lebesgue measure, satisfying the following estimate: for each $t_0>0$, there exists  a constant $C(t_0)$ independent of the initial value $X_0=x$ such that
$$\sup_{t\geq t_0}\|p_t\|_{B_{1,\infty}^k}\leq C(t_0)<\infty,\quad \text{for each } t_0>0.$$
\end{corollary}

\begin{proof}
Having checked the Malliavin differentiability of $X_t$ up to order $k+1$ and invertibility of its Malliavin matrix, existence and regularity of $p_t$ follow by the previous results. The temporal estimate for $t\geq t_0$ then follows from a combination of Claim \ref{claims0}, the second half of Proposition \ref{claims1} and Claim \ref{claims2}.
\end{proof}

In this paper Malliavin calculus is employed as a tool to establish regularity estimates of the density. While one usually requires the SDE to start from a fixed point $X_0=x$ to establish Malliavin differentiability, it is customary in McKean-Vlasov dynamics to consider general initial distributions that are not atomic. The next corollary settles this issue.

\begin{corollary}\label{corollary3.13meaning}
Consider the SDE $dX_t=a(t,X_t)dt+dW_t$ with initial distribution $\mu_0\in\mathcal{P}(\mathbb{R}^d).$ Assume that $a\in L^\infty([0,T]\times\mathbb{R}^d).$

\begin{enumerate}
\item Assume for some $k\in\mathbb{N}_+$, $Da,D^2a,\cdots,D^ka\in L^\infty([0,T]\times\mathbb{R}^d)$. Then the same conclusion as in Corollary \ref{corollarystep1} hold true, i.e. for each $t>0$, $X_t$ has a density $p_t$ with respect to the Lebesgue measure, satisfying
\begin{equation}
    \label{what'snew?}
\sup_{t\geq t_0}\|p_t\|_{B_{1,1}^k}\leq C(t_0)<\infty,\quad \text{for each } t_0>0.\end{equation}
\item More generally, assume that for each $t_0>0$, $Da,\cdots,D^ka\in L^\infty([t_0,T]\times\mathbb{R}^d)$, that is, the gradients are bounded after any positive time. Then the same conclusion holds.
\end{enumerate}
\end{corollary}

\begin{proof}
For each $x\in\mathbb{R}^d$ and $t>0$, denote by $p^x_t$ the law of the SDE $dX_t=a(t,X_t)dt+dW_t$  at time $t$, with initial value $X_0=x$. Then by Corollary \ref{corollarystep1}, each $p_t^x$ satisfies \eqref{what'snew?} with $p_t^x$ in place of $p_t$, having the same constant $C(t_0)$. By Markov property of the SDE $X_t$ we have
$$p_t(y)=\int_{\mathbb{R}^d}p_t^x(y) \mu_0(dx),$$
then by Jensen's inequality, 
$$\|p_t\|_{B_{1,1}^k}\leq \int_{\mathbb{R}^d}\|p_t^x\|_{B_{1,1}^k}\mu_0(dx)<C(t_0)<\infty.$$

This proves the first claim. To prove the second claim,
for each $t_0>0$ we fix some $0<s_0<t_0$. By Markov property, the law of $X_{t_0}$ is the same as the law of $Y_{t_0-s_0}$, where $Y$ solves the SDE $dY_t=a(t-s_0,Y_t)dt+dW_t$ with initial distribution $\operatorname{Law}(Y_0)=\operatorname{Law}(X_{s_0}).$ Then the second claim follows from the first.
\end{proof}

\subsection{Applications to McKean-Vlasov SDE}\label{section3.5}
We have made all technical preparations for the proof of Theorem \ref{theorem2}. 

\begin{theorem}
Consider the McKean-Vlasov SDE 
\begin{equation}\label{theorem3.14}
    dX_t=\langle b(t,X_t-\cdot),\mu_t\rangle dt +dW_t,\quad \operatorname{Law}(X_t)=\mu_t
\end{equation} with initial law $\mu_0\in\mathcal{P}(\mathbb{R}^d).$ Assume we have $b(t,\cdot)\in \mathcal{C}_b^1(\mathbb{R}^d)$ for each $t>0$ with the estimate
$$\sup_{t\geq 0}\|b(t,\cdot)\|_{\mathcal{C}_b^1}<\infty.$$Then for each $t>0$ the law of $X_t$ has a \textbf{smooth} density $p_t$ with respect to the Lebesgue measure.
\end{theorem}

\begin{proof}
Fix an arbitrarily large $T_0>0$.
We proceed by induction and prove the following claim for each $k\in\mathbb{N}_+$: for each $t_0>0$ there exists a constant $C(t_0)<\infty$ such that  
\begin{equation}\label{inductionstep1}
  \sup_{t_0\leq t\leq T_0}  \|p_t\|_{B_{1,1}^{k}}<C(t_0)<\infty\quad \text{ for all } t_0>0.\end{equation}
  
  Since $b\in \mathcal{C}_b^1$, 
  the function $x\mapsto \langle b(t,x-\cdot),\mu\rangle$ is continuously differentiable with bounded gradient for any probability measure $\mu\in\mathcal{P}(\mathbb{R}^d)$. Then we can apply Corollary \ref{corollary3.13meaning} with $k=1$ and justify \eqref{inductionstep1} with $k=1$. Suppose that \eqref{inductionstep1} is justified with exponent $k=\ell\in\mathbb{N}^+$, we have the following estimate, where we regard each function as a function in $x$ when taking the norm:
  $$
  \begin{aligned}
  \|\langle b(t,x-\cdot,\mu_t\rangle\|_{B_{\infty,1}^{\ell+1}}&=\left\|\int_{\mathbb{R}^d}b(t,x-y)p_t(y)dy\right\|_{B_{\infty,1}^{\ell+1}}\\ &\leq \|b(t,\cdot)\|_{B_{\infty,\infty}^{1}}\|p_t\|_{B_{1,1}^\ell}\\&\leq \|b(t,\cdot)\|_{\mathcal{C}_b^{1}}C(t_0)<\infty\quad \text{ for all }t>t_0.
  \end{aligned}
  $$
 In the second line we used Young's convolution inequality \eqref{convolutioneq2}. By \eqref{functionrelations}, the first $\ell+1$-th derivatives of the drift of the McKean-Vlasov SDE \eqref{convolutioneq2} lie in $B_{\infty,1}^0\hookrightarrow L^\infty$, thus verifying assumptions in the second part of Corollary \ref{corollary3.13meaning}, with $k=\ell+1$. By that corollary, we justify \eqref{inductionstep1} with exponent $k=\ell+1.$

By induction, we have shown that $p_t\in\cap_{k\in\mathbb{N}_+}B_{1,1}^k\hookrightarrow \cap_{k\in\mathbb{N}_+}B_{1,\infty}^k$. By the Sobolev embedding $B_{1,\infty}^k\hookrightarrow B_{\infty,\infty}^{k-d}\hookrightarrow \mathcal{C}_b^{k-d-\gamma}$ for any $\gamma>0$, we deduce that $p_t$ is smooth.
\end{proof}

\begin{remark}\label{remark2}
The assumption $b\in\mathcal{C}_b^{1}$ comes from the fact that integration by parts via Malliavin calculus does not increase regularity: an SDE with $\mathcal{C}^k$ drifts yield at most $k$ time smoothness in Besov sense (see Corollary \ref{corollary3.13meaning}), so we have to take extra regularity into the system otherwise our argument doesn't work. In a first version of this draft, we made a stronger assumption that $b\in \mathcal{C}_b^{1+\alpha}$ for some $\alpha>0$, in the hope that we can avoid encountering the spaces $B_{\infty,\infty}^s$, $s\in\mathbb{N}_+$ which contain functions that are not $\mathcal{C}^k$. Indeed, the results in Section \ref{section3.2} and \ref{section3.3} only require that the coefficients lie in the Sobolev space $W^{k,\infty}$, which is a weaker condition than $\mathcal{C}^k$. Combined with a better use of inequality \eqref{convolutioneq2}, the $\alpha>0$ assumption seems unnecessary.
\end{remark}

\section{Riesz transform: towards a better regularity}\label{section4}

In this section we turn back to the $L^\infty$ case and improve the regularity result we just obtained, completing the proof of Theorem \ref{theorem1}. We also prove Corollary \ref{corollary1.3} in Section \ref{section4.3}.

\subsection{A review of the framework} 

We first give a brief review of the framework based on Riesz transform, introduced in \cite{bally2011riesz}.

Recall that the fundamental solution of the Laplace equation $\Delta Q_d=\delta_0$ in $\mathbb{R}^d$ is
$$Q_2(x)=a_2^{-1}\ln x,\quad Q_d(x)=a_d^{-1}|x|^{2-d},d>2.$$ 
Then for $f\in \mathcal{C}_c^1(\mathbb{R}^d)$ we have the following representation, named Riesz transform:
$$f=(\nabla Q_d)*\nabla f.$$

Following \cite{bally2011riesz}, for a Borel probability measure $\mu$ on $\mathbb{R}^d$ we define $$L_\mu^p:=\{\phi:\int|\phi(x)|^p\mu(dx)<\infty\},$$ and for $\phi \in L_\mu^p$, define $\partial_i^{\mu}\phi$, if exists, as the unique $\theta_i$ satisfying 
$$\int \partial_i f(x)\phi(x)\mu(dx)=-\int f(x)\theta_i(x)\mu(dx),\quad f\in\mathcal{C}_c^\infty(\mathbb{R}^d).$$
From this we define the Sbolev space $W_\mu^{1,p}$ as the completion under the norm $$\|\phi\|_{W_\mu^{1,p}}=\|\phi\|_{L_\mu^p}+\sum_{i=1}^d \|\partial_i^\mu\phi\|_{L_\mu^p}.$$ 
Higher order Sobolev spaces $W_\mu^{m,p}$, $m\geq 2$ can be defined analogously. We have the following useful representation formula of probability density functions (see \cite{bally2011riesz}, Theorem 1): 

Assume that $1\in W_\mu^{1,p}, p\geq 1$, then $\mu(dx)=p_\mu(x)dx$ where
$$p_\mu(x)=-\sum_{i=1}^d \int \partial_i Q_d(y-x)\partial_i^\mu 1(y)\mu(dy),$$ and we have 
$\partial_i p_\mu(x)=p_\mu(x)\partial_i^\mu 1(x)$.

We will use the following crucial estimate of Riesz transform:

define $$\Theta_p(\mu)=\sup_{a\in\mathbb{R}^d}\sum_{i=1}^d \left(\int_{\mathbb{R}^d}|\partial_i Q_d(x-a)|^{\frac{p}{p-1}}\mu(dx)\right)^{\frac{p-1}{p}}.$$
Then by Theorem 5 and Theorem 8 of \cite{bally2011riesz}, for $p>d$ and $m$ with $1\in W_\mu^{m,p}$, we have 
\begin{equation}\Theta_p(\mu)\leq dK_{dp}\|1\|_{W_\mu^{1,p}}^{k_{d,p}},\quad
\|p_\mu\|_\infty\leq 2dK_{dp}\|1\|_{W_\mu^{1,p}}^{k_{d,p}+1},\end{equation}
\begin{equation}
\label{higherorderestimate}    
\|p_\mu\|_{W^{m,p}}\leq (2dK_{d,p})^{1-1/p}\|1\|_{W_\mu^{1,p}}^{k_{d,p}(1-1/p)}\|1\|_{W_\mu^{m,p}}.\end{equation}
We refer to \cite{bally2011riesz} for values of the numeric constants $K_{dp}$ and $k_{dp}.$ 

Now turn back to the setting of Malliavin calculus. Consider a random variable $X$ satisfying Lemma \ref{lemma3.6}, that is, $X\in\mathcal{S}^k$ and $\gamma_X^{-1}\in\mathcal{S}^k$ for some $k\geq 2$, and consequently there exists $\bar{Z}\in \cap_{p\geq 1}L^p(\Omega)$ such that $\mathbb{E}[\partial_i G(X)]=\mathbb{E}[G(X)\bar{Z}]$ for any $G\in\mathcal{C}_b^\infty.$ If in the aforementioned framework we let $\mu$ 
 be the law of $X$, then we have $\partial_i^\mu 1=\bar{Z}\in \cap_{p\geq 1} L^p(\Omega)$, so that $1\in \cap_{p\geq 1} W_\mu^{1,p}$. This implies the density function $p_\mu\in W^{1,p}$ for any $p\geq 1$,  with $\|p_\mu\|_{W^{1,p}}$ bounded by $\|\bar{Z}\|_{L_\mu^p}$. 
 
Considering the Sobolev spaces $W^{1,p}$, $p>d$ has the benefit that we can prove Hölder regularity of $p_\mu$ which is not attainable for $W^{1,p}$ with smaller values of $p$, as suggested in section 2.4 of \cite{bally2011riesz}. Indeed, by Morrey's inequality, for $p>d$ we have continuous embeddings $W^{1,p}\hookrightarrow \mathcal{C}_b^{1-\frac{d}{p}}$, and more generally $W^{m,p}\hookrightarrow \mathcal{C}_b^{m-\frac{d}{p}}$ for $m\geq 1$. Therefore if $p_\mu\in W^{1,p}$ for all $p\geq 1$, we deduce that $p_\mu$ is $\alpha$-Hölder for any $\alpha\in(0,1).$
 
We make a final remark that by Meyer inequality \eqref{meyerinequality} and Hölder inequality, $\|\bar{Z}\|_{L_\mu^p}$ is further bounded by $\|X\|_{\mathbb{D}^{1,p'}}$ and $\|\gamma_X^{-1}\|_{\mathbb{D}^{1,p''}}$ for some $p',p''>p.$

\subsection{Application to McKean-Vlasov SDE}\label{section4.2}
Consider the McKean-Vlasov SDE $$dX_t=\langle b(t,X_t-\cdot),\mu_t\rangle dt+dW_t\quad\operatorname{Law}(X_t)=\mu_t$$ with $b\in L^\infty.$ We have shown in Corollary \ref{finalcorollary} that for each $t>0$, $X_t$ admits a density $p_t$ with respect to the Lebesgue measure, such that for each $\alpha\in(0,1)$ and each $t_0>0$, we can find a constant $C(t_0)$ such that
$$\sup_{t\geq t_0} \|p_t\|_{B_{1,\infty}^{2+\alpha}}<C(t_0)<\infty,\quad t_0>0.$$

In dimension 1, the Sobolev embedding $B_{1,\infty}^s\hookrightarrow B_{\infty,\infty}^{s-1}$ implies $p_t\in\mathcal{C}_b^{1+\alpha}$ for each $\alpha\in(0,1)$. Now we show that $p_t\in\mathcal{C}_b^{1+\alpha}$  holds in any dimension. 

Set $a(t,x):=\langle b(t,x-\cdot,\mu_t\rangle.$ The estimate on $p_t$ implies 
$$\sup_{t\geq t_0} \|a(t,\cdot)\|_{B_{\infty,\infty}^{2+\alpha}}<C(t_0)<\infty,\quad t_0>0.$$

For any $t>0$ fix some $t_0\in(0,t).$
Consider the SDE $dY_t=a(t-t_0,Y_t)dt+dW_t$ with initial value $Y_0=x\in\mathbb{R}^d$. By  results of Malliavin calculus (Proposition \ref{propositionuseful} and \ref{claims1}), $Y_{t-t_0}$ is (three times) Malliavin differentiable, whose Malliavin matrix has inverse moments of all orders. Denote by $p_{t,t_0}^x$ the law of $Y_{t-t_0}$ and set $\mu:=p_{t,t_0}^x$. Observe that by Lemma \ref{lemma3.6} with $k=3$, the variable $\bar{Z}$ (which is equal to  $\partial_i^\mu 1$ by definition) lies in $\mathcal{S}^1$. By applying Lemma \ref{lemma3.6} the second time we see that $1\in W_\mu^{2,p}$ for each $p\geq 1$.

Now we can use equation \eqref{higherorderestimate} with $m=2$ to deduce that $p_{t,t_0}^x\in W^{2,p}$ for every $p>d$, and  $\|p_{t,t_0}^x\|_{W^{2,p}}$ can be bounded uniformly in $x$. By the Markov property of $X_t$, we have $$p_t(y)=\int_{\mathbb{R}^d} p_{t,t_0}^x(y) p_{t_0}(x)dx,\quad y\in\mathbb{R}^d,$$ then by Jensen's inequality we have $p_t\in W^{2,p}$ for each $p>d$ and every $t>0$. Consequently $p_t\in\mathcal{C}_b^{1+\alpha}$ for each $\alpha\in(0,1)$ and $t>0$ by Morrey's inequality. As claimed, the regularity of $p_t$ is strictly better than the example in in Appendix \ref{appendixB}, in any dimension.

This completes the proof of Theorem \ref{theorem1}.

\subsection{For Hölder interactions}\label{section4.3}

In this subsection we prove Corollary \ref{corollary1.3}.

\begin{proof}
By the result of Corollary \ref{finalcorollary}, $p_t$, the density function of $\mu_t$, lies in $B_{1,\infty}^{2+\gamma}$ for each $\gamma>0$. Since $b\in\mathcal{C}_b^\alpha$ for some $\alpha>0$, we can find some $\theta>0$ such that 
$$\sup_{s\geq t_0>0}\|x\mapsto \langle b(s,x-\cdot),\mu_s\rangle\|_{\mathcal{C}_b^{3+\theta}}<C(t_0)<\infty,$$
that is, the drift of the McKean-Vlasov SDE is $\mathcal{C}^3$ in space. Now carrying out the same argument as in Section \ref{section4.2}, but with $m=3$ instead of $m=2$, we can show that $p_t\in\mathcal{C}_b^{2+\gamma}$ for any $\gamma\in(0,1)$ and any $t>0$.
\end{proof}

\section{More examples}\label{section5}
In this section we show that Theorem \ref{theorem1} and \ref{theorem2} hold under more general assumptions on the kernel $b$ and more general driving noise $Z$ other than Brownian motion.
\subsection{Singular interaction kernels}\label{section5.1}
One extension is to consider kernels $b$ that are not locally bounded. Consider $b:[0,T]\times\mathbb{R}^d$ that satisfies the integrability condition \begin{equation}\label{5.0}b\in L^q([0,T];L^p(\mathbb{R}^d)),\quad p\geq 1,q\geq 2, \frac{d}{p}+\frac{2}{q}<1.\end{equation}
Well-posedness of the SDE 
\begin{equation}\label{5.1}dX_t=b(t,X_t)dt+dW_t\end{equation} under \eqref{5.0} is established in \cite{krylov2005strong}, 
and well-posedness of the McKean-Vlasov SDE, with the same assumption on $b$,
\begin{equation}\label{5.2}dX_t=\langle b(t,X_t-\cdot,\mu_t\rangle dt+dW_t,\quad \mu_t=\operatorname{Law}(X_t)\end{equation}
is established in \cite{rockner2021well}, see also \cite{han2022solving}.

For any collection of probability measures $(\mu_t)_{t\in[0,T]}$ with each $\mu_t\in\mathcal{P}(\mathbb{R}^d)$, one observes by Jensen's inequality that the function
$$(t,x)\mapsto \langle b(t,x-\cdot),\mu_t\rangle \in L^q([0,T];L^p(\mathbb{R}^d)).$$ Therefore we may regard \eqref{5.2} as a special case of \eqref{5.1}, once the measure $\mu_t$ is known.

Now we can apply Theorem 5.1 of \cite{romito2018simple}, which states that

\begin{proposition}
\label{proposition5.1}
Assume $b$ satisfies \eqref{5.0}, consider the SDE \eqref{5.1} with $X_0=x\in\mathbb{R}^d$, then for each $t>0$ $X_t$ has a density $p_t$ with respect to the Lebesgue measure, such that for any $\gamma<1-\frac{2}{q}$ and $e_\gamma>\frac{1-\frac{1}{q}}{1-\frac{2}{p}-\frac{d}{q}}\gamma$, 
$$\sup_{t\in[0,T]}(1\wedge t)^{e_\gamma}\|p_t\|_{B_{1,\infty}^\gamma}<\infty.$$
\end{proposition}

Since $b_t\in L^p\hookrightarrow B_{p,\infty}^0$, we conclude that for each $t$, $\langle b(t,x-\cdot,\mu_t\rangle\in B_{p,\infty}^\gamma(\mathbb{R}^d)$ as a function of $x$, via the convolution inequality \eqref{convolutioninequ}.

By Sobolev embedding, $B_{p,\infty}^s\hookrightarrow B_{\infty,\infty}^{s-\frac{d}{p}}$ for any $s\in\mathbb{R}$. Since $\frac{d}{p}+\frac{2}{q}<1$, we can choose $\gamma$ sufficiently close to $1-\frac{2}{q}$ and obtain that $\langle b(t,x-\cdot,\mu_t\rangle\in B_{\infty,\infty}^\delta(\mathbb{R}^d)$ for any $\delta\in(0,1-\frac{d}{p}-\frac{2}{q})$. That is, the drift of the McKean-Vlasov SDE is $\delta$-Hölder continuous in space. 
\begin{corollary}
Assume $b$ satisfies \eqref{5.0} and consider the McKean-Vlasov SDE \eqref{5.2}. Then for each $t>0$, the function $x\mapsto \langle b(t,x-\cdot),\mu_t\rangle$ is $\delta$-Hölder continuous on $\mathbb{R}^d$, for each $\delta\in(0,1-\frac{d}{p}-\frac{2}{q}).$
\end{corollary}

In this subsection we assume $q=\infty$, that is, 
\begin{equation}\label{5.4}
b\in L^\infty([0,T],L^p(\mathbb{R}^d)),\quad p>d.\end{equation}
This simplifies the computations and allows us to apply results proved in Section \ref{section2}. The case of general $q>2$ will be dealt with in Section \ref{section5.20}.

Now by what we have shown, \begin{equation}
    \label{first12}
\sup_{t\geq t_0}\|x\mapsto \langle b(t,x-\cdot),\mu_t\rangle\|_{B_{\infty,\infty}^\delta}\leq C(t_0)<\infty,\quad t_0>0,\end{equation}
so we can apply results in Section \ref{section2.3} and deduce that 
$$\sup_{t\geq t_0}\|p_t\|_{B_{1,\infty}^{1+\gamma}}\leq C(t_0,\gamma)<\infty,\quad t_0>0,$$
for each $0<\gamma<\delta.$

This implies a bound on $\|x\mapsto \langle b(t,x-\cdot),\mu_t\rangle\|_{B_{p,\infty}^{1+\gamma}}$ since $b\in L^p$. By Sobolev embedding $B_{p,\infty}^s\hookrightarrow B_{\infty,\infty}^{s-\frac{d}{p}}$, we deduce the following bound:
\begin{equation}\label{second12}
\sup_{t\geq t_0}\|x\mapsto \langle b(t,x-\cdot),\mu_t\rangle\|_{B_{\infty,\infty}^{1+\gamma-\frac{d}{p}}}\leq C(t_0)<\infty,\quad t_0>0.\end{equation}
Compared with \eqref{first12}, since $1-\frac{d}{p}>0$, we have improved the regularity parameter of the drift by a fixed positive value. 

We carry out the procedure a number of times, until $x\mapsto \langle b(t,x-\cdot),\mu_t\rangle$ is in $\mathcal{C}_b^\theta$ for some $\theta>1$. Then we apply results in Section \ref{section2.4} a number of times, until $x\mapsto \langle b(t,x-\cdot),\mu_t\rangle$ is in $\mathcal{C}_b^\theta$ for some $\theta>2$. 

The last step, i.e. showing $p_t\in \mathcal{C}_b^{1+\alpha}$, is identical to the proof of the second part of Theorem \ref{theorem1}. We turn to Section \ref{section4}, use Malliavin calculus (because the drift is now $C^2$, $X_t$ is three times Malliavin differentiable by Proposition \ref{propositionuseful}), and the estimates based on Riesz transform (as we did in Section \ref{section4.2}, we show $1\in W_\mu^{2,p}$, $p\geq 1$, then apply \eqref{higherorderestimate} with $m=2$, and Morrey's inequality. Alternatively, we can apply the main result of \cite{banos2016malliavin} directly), to finally deduce that $p_t\in\mathcal{C}_b^{1+\alpha}$ for each $\alpha\in(0,1).$

This proves Theorem \ref{theorem3} in the special case $q=\infty.$

\subsection{Time dependent singular drifts}\label{section5.20}
We begin with the following elementary estimate:
\begin{lemma}\label{lemma5.3}
Consider the SDE 
$$dX_t=a(t,X_t)dt+dW_t$$
where $a:[0,T]\times\mathbb{R}^d\to\mathbb{R}^d$ satisfies $a(t,\cdot)\in L^\infty$ for each $t\in[0,T]$ whose norm further satisfies the following integrability condition: for some $q>2$, one has $t\mapsto \|a(t,\cdot)\|_\infty\in L^q([t_0,T])$  for each $t_0>0$. Then we have the estimate 
$$\mathbb{E}[|X_s-X_t|]\leq C(t_0)|s-t|^{\frac{1}{2}},\quad |s-t|\leq 1, s,t\geq t_0>0.$$
\end{lemma}

\begin{proof}
It is easy to see that $\mathbb{E}[|W_s-W_t|]\leq C_d |s-t|^{\frac{1}{2}}$ for some $C_d$ depending only on $d$. Moreover, 
$$\left|\int_s^t a(r,X_r)dr\right|\leq |t-s|^{\frac{1}{\tilde{p}}}\left(\int_s^t \|a(r,\cdot)\|_{L^\infty}^{q} dr\right)^\frac{1}{q}<C(t_0)|t-s|^{\frac{1}{2}},$$
where $\tilde{p}=\frac{q}{q-1}\in(1,2)$ is the conjugate exponent of $q$.  
\end{proof}

Now we extend the estimates in Section \ref{section2.3} and \ref{section2.4} to drifts of SDE that depend on time via an integrability condition. As before, we work around with the baby model (where $a$ satisfies the assumptions given in Lemma \ref{lemma5.3})
$$dX_t=a(t,X_t)dt+dW_t.$$

For the zeroth order approximation, we begin with 
$$Y_s=X_{t-\epsilon}+W_s-W_{t-\epsilon},\quad s\geq t-\epsilon.$$

The approximation error, \textbf{AE}, satisfies (for $\alpha\in(0,1)$),
$$\textbf{AE}\leq[\phi]_{\mathcal{C}_b^\alpha}\mathbb{E}\left[\left|\int_{t-\epsilon}^t a(s,X_s)ds\right|^\alpha\right]\leq C(t_0) [\phi]_{\mathcal{C}_b^\alpha}\epsilon^{\frac{\alpha}{\tilde{p}}},$$
where $\tilde{p}$ is the conjugate exponent of $q$.

The probability estimate, \textbf{PE}, satisfies 
$$\textbf{PE}:=\mathbb{E}[\Delta_h^m \phi(Y_t^\epsilon)]\leq \|\phi\|_{L^\infty} (|h|/\sqrt{\epsilon})^m,$$

Gathering the approximations \textbf{AE} and \textbf{PE}, utilizing Proposition \ref{technicalproposition}, we see that $X_t$ has a law $p_t$ that lies in $B_{1,\infty}^\gamma$ for all $\gamma\in(0,1-\frac{2}{q})$. This is consistent with Proposition \ref{proposition5.1}.  

For the first order approximation, we assume that for some $\beta\in(0,1)$, $a(t,\cdot)\in \mathcal{C}_b^\beta$, and the norm satisfies the temporal estimate $t\mapsto \|a(t,\cdot)\|_{\mathcal{C}_b^\beta}\in L^q([t_0,T])$ for any $t_0>0$. Defining $Y_t$ as in Section \ref{section2.3}:
$$Y_s=X_{t-\epsilon}+\int_{t-\epsilon}^s a(r,X_{t-\epsilon})dr+W_s-W_{t-\epsilon},\quad s\geq t-\epsilon,$$
then we have 
$$
\begin{aligned}
\mathbb{E}\left[|X_t-Y_t|\right]=&\mathbb{E}\left|\int_{t-\epsilon}^t (a(r,X_x)-a(r,X_{t-\epsilon}))dr\right|\\ &\leq \int_{t-\epsilon}^t  \|a(r,\cdot)\|_{\mathcal{C}_b^\beta}\mathbb{E}[\|X_r-X_{t-\epsilon}\|^\beta] ds\leq \epsilon^{\frac{\beta}{2}+\frac{1}{\tilde{p}}},\end{aligned}$$
where in the last step we used Lemma \ref{lemma5.3}. Thus the approximation error, \textbf{AE}, is given by 
$$\textbf{AE}\leq [\phi]_{\mathcal{C}_b^\alpha}\epsilon^{\alpha(\frac{\beta}{2}+\frac{1}{\tilde{p}})}.$$
The probabilistic estimate, \textbf{PE}, is the same as in Section \ref{section2.3}. Then using Proposition \ref{technicalproposition}, we see that the density function $p_t$ of $X_t$ lies in $B_{1,\infty}^\gamma$ for any $\gamma\in(0,\beta+1-\frac{2}{q})$.

For the second order approximation, we assume that for some $\beta\in(0,1)$, $a(t,\cdot)\in\mathcal{C}_b^{1+\beta}$ with the temporal estimate $t\mapsto \|a(t,\cdot)\|_{\mathcal{C}_b^{1+\beta}}\in L^q([t_0,T])$ for any $t_0>0$. Defining $Y_t$ as in Section \ref{section2.4}:
$$Y_s=X_{t-\epsilon}+\int_{t-\epsilon}^s \left[a(r,X_{t-\epsilon})+Da(r,X_{t-\epsilon})(W_r-W_{t-\epsilon})\right] dr
+W_s-W_{t-\epsilon},\quad s\geq t-\epsilon.$$
Now computing $Y_s-X_s$, it is easy to check that 
$$\mathbb{E}[|Y_t-X_t|]\leq\int_0^\epsilon \|a(t-\epsilon+r,\cdot)\|_{\mathcal{C}_b^{1+\beta}} r^{\min(\frac{1+\beta}{2},\frac{1}{\tilde{p}})} dr\leq C \epsilon ^{\min(\frac{1+\beta}{2},\frac{1}{\tilde{p}})+\frac{1}{\tilde{p}}},$$
where again $\tilde{p}$ is the conjugate exponent of $q$. The probability estimate, \textbf{PE}, is the same as in Section \ref{section2.4}. Then applying Proposition \ref{technicalproposition}, we see that $p_t$ lies in $B_{1,\infty}^\gamma$ for each $\gamma\in (0,\min(\beta+2-\frac{2}{q},3-\frac{4}{q})).$ 

At this point, one sees that the value of $q$ significantly affects the potential of our bootstrap argument. Say if $q\in(2,4]$, we must have $\gamma<2$, and the bootstrap argument cannot proceed anymore. While if $q=\infty$, we can choose $\gamma>2$ freely, and that is the case of Section \ref{section2.4}.

\begin{corollary}
For the McKean-Vlasov SDE \eqref{1.43}, assuming $b\in L^q([0,T];L^p(\mathbb{R}^d))$ for some $\frac{d}{p}+\frac{2}{q}<1$, $p>d$, $q>2$. Then for each $t>0$, $\mu_t$ has a density $p_t$ that lies in $B_{1,\infty}^{\gamma}$ for any $\gamma\in(0,3-\frac{4}{q})$.
\end{corollary}
\begin{proof} One just needs to apply Proposition \ref{proposition5.1} and the three estimates in this section, to improve regularity of $p_t$ via a bootstrap argument until it reaches $3-\frac{4}{q}$. Indeed, when the drift is known to be $\mathcal{C}_b^\gamma$, one proves that $p_t\in B_{1,\infty}^{1+\gamma-\frac{2}{q}-}$, by the various estimates in this chapter. Then by convolution the drift is in $B_{p,\infty}^{1+\gamma-\frac{2}{q}-}\hookrightarrow B_{\infty,\infty}^{1+\gamma-\frac{2}{q}-\frac{d}{p}-}$, so we are improving the regularity by a positive factor each time since $1>\frac{2}{q}+\frac{d}{p}$.
\end{proof}
This completes the proof of Theorem \ref{theorem3}.

\subsection{Distributional interaction kernels}\label{section5.30}

Now consider the McKean-Vlasov SDE \eqref{1.43} where the interaction $b\in \mathcal{C}_b^\alpha:=B_{\infty,\infty}^\alpha$ for some $\alpha\in(-1,0).$ The well-posedness issue has been addressed recently in \cite{de2022multidimensional}. We present a particular, yet sufficiently general case established in \cite{de2022multidimensional}, whose main result can be stated as:

\begin{theorem}
Assume $b\in L^\infty([0,T];B_{\infty,\infty}^\alpha(\mathbb{R}^d))$ for $\alpha\in(-1,0)$. Then for any $\mu_0\in\mathcal{P}(\mathbb{R}^d)$, the McKean-Vlasov SDE \eqref{1.43} has a unique strong solution, satisfying that for each $t>0$, the law $\mu_t$ has a density $p_t$ with respect to the Lebesgue measure, $p_t\in B_{1,1}^{-\alpha}$, and the temporal estimate $p_t\in L^{r}([0,T];B_{1,1}^{-\alpha})$
holds for each $r\in[1,-\frac{2}{\alpha}).$
\end{theorem}

By the convolution inequality \eqref{convolutioneq2}, the drift of the McKean-Vlasov SDE $\langle b(t,x-\cdot),\mu_t\rangle\in B_{\infty,1}^0(\mathbb{R}^d)\hookrightarrow L^\infty(\mathbb{R}^d)$, with the temporal estimate $\langle b(t,x-\cdot),\mu_t\rangle\in L^r([0,T];L^\infty(\mathbb{R}^d))$ for each $r\in [1,-\frac{2}{\alpha})$. By choosing $r$ sufficiently close to $-\frac{2}{\alpha}$, one sees that the McKean-Vlasov SDE \eqref{1.43} can be regarded as a SDE with drift satisfying \eqref{5.0}, so Proposition \ref{proposition5.1} is applicable. 

By Proposition \ref{proposition5.1}, $p_t\in B_{1,1}^{\gamma}$ for each $\gamma\in (0,1-\frac{2}{-\frac{2}{\alpha}})=(0,1+\alpha)$ with the temporal estimate $\|p_t\|_{B_{1,1}^\gamma}\in L^\infty([t_0,T])$ for each $t_0>0>0$ and $T<\infty$.

Then by the convolution inequality \eqref{convolutioneq2}, the drift $\langle b(t,x-\cdot),\mu_t\rangle\in B_{\infty,1}^{\gamma+\alpha}(\mathbb{R}^d)$. By choosing $\gamma$ sufficiently close to $1+\alpha$, one sees that $\alpha=-\frac{1}{2}$ is a critical threshold, below which the drift of the McKean-Vlasov SDE lies in negative index Besov spaces and our bootstrap fails.

\begin{itemize}
\item $\alpha\in(-\frac{1}{2},0]$, then the regularity of $p_t$ is improving. Using the estimates in Section \ref{section5.20}, we can eventually show $p_t\in B_{1,\infty}^{\gamma}$ for each $\gamma\in (0,3)$. Indeed, we first show the drift of the McKean-Vlasov SDE is Hölder continuous, with Hölder norm lying in $L^\infty([t_0,T])$, for any $t_0>0.$ Then each step via the bootstrap argument (see Section \ref{section2.3}), we improve the regularity of $p_t$ by a factor that is arbitrarily close to $1$, and by convolution with a drift $b\in \mathcal{C}_b^{\alpha}$, the regularity of the drift of the SDE is improved by a factor arbitrarily close to $\alpha+1>0$.
Apply this procedure a number of times, until the drift lies in $\mathcal{C}_b^{1+\gamma}$ for some $\gamma>0$, and we get along the way that its $\mathcal{C}_b^{1+\gamma}$ norm is in $L^\infty([t_0,T])$ for any $t_0>0$. Now proceed as in Section \ref{section2.4}, we improve the regularity of the drift until it is more than $\mathcal{C}^2$. Along the way we have shown that $p_t\in B_{\infty,\infty}^\gamma$ for $\gamma\in(0,3)$, and using the argument of Section \ref{section4} we show that $p_t$ is also in $\mathcal{C}_b^{\gamma}$ for any $\gamma<2$. This proves Theorem \ref{theorem1.7}.

\item $\alpha\in (-1,-\frac{1}{2}]$. This is the range where our bootstrap argument fails, and also (surprisingly) the range where well-posedness of the SDE $dX_t=V(X_t)dt+dW_t$ may fail. The classical explanation for the $-\frac{1}{2}$ threshold is that this is the place where the paraproduct structure in the corresponding backward Kolmogorov equation may be ill defined, see for example \cite{flandoli2017multidimensional} and the introduction of \cite{de2022multidimensional}. We do however mention that if $V$ has an additional rough path structure, the range of solvability can be extended to $V\in B_{\infty,\infty}^\gamma$ for $\gamma>-\frac{2}{3}$, see for example \cite{cannizzaro2018multidimensional} and \cite{kremp2022multidimensional}.
\end{itemize}

\subsection{Stable driven SDEs}\label{section5.2}

Let $(Z_t)_{t\geq 0}$ be an $\alpha$-stable-like process with $\alpha\in(1,2)$. It is proved in \cite{debussche2013existence} and Section 7.2 of \cite{romito2018simple} that with bounded measurable drift and a nondegenerate stable noise, the density exists and lies in some Besov spaces. More precisely,

\begin{theorem}[See \cite{debussche2013existence} or Theorem 7.5 of \cite{romito2018simple}] Given $(Z_t)_{t\geq 0}$ an $\alpha$-stable process with $\alpha\in(1,2)$, and $a\in L^\infty([0,T]\times\mathbb{R}^d),$ consider the SDE $$dX_t=a(t,X_t)dt+dZ_t$$ starting from $x$. Then for each time $t>0$ the SDE has a density $p_x(t)$ with respect to the Lebesgue measure of $\mathbb{R}^d$, satisfying the following estimate: for any $\gamma\in(0,\alpha-1)$,
$$\|p_x(t)\|_{B_{1,\infty}^\gamma}<\infty.$$
\end{theorem}

Fix a kernel $b\in L^\infty([0,T]\times\mathbb{R}^d)$, consider the stable driven McKean-Vlasov SDE
\begin{equation}\label{stablemckeanvlasov}
    dX_t=\langle b(t,X_t-\cdot),\mu_t\rangle dt+dZ_t,\quad \operatorname{Law}(X_t)=\mu_t
\end{equation} 
that starts from $x\in\mathbb{R}^d$. Some very general well-posedness results of stable driven McKean-Vlasov SDEs can be found in \cite{de2022multidimensional}.

By the previous theorem, we conclude that for each $t>0$, $\mu_t\in B_{1,\infty}^{\gamma}$ for any $\gamma<\alpha-1$. Therefore the function $x\mapsto \langle b(t,x-\cdot),\mu_t\rangle$ is also $\gamma$-Hölder continuous for $\gamma<\alpha-1$. Since $\alpha>1$, we have improved regularity of the drift, and this procedure can possibly be carried over a number of times to show $p_x(t)\in\mathcal{C}_b^{1+\alpha}$, $\alpha>0$.

Once we have set up the corresponding Malliavin calculus toolbox for stable processes, that is comparable to the ones in Section \ref{section2} for Brownian motion, an analogous result of Theorem \ref{theorem2} can possibly be proved.

\appendix
\section{Technical proofs}\label{appendixA}
We now give a proof of Proposition \ref{technicalproposition}.

\begin{proof}
Given $a<a_0$, we first choose some $\delta_1<a_0-a$ such that $\alpha:=\frac{a+\delta_1}a_0<1.$ This is the value of $\alpha$ we use in the proof. Find another $\delta_2\in(0,1)$ that will be defined later, and choose $m$ sufficiently large such that 
\begin{equation}\label{assumptiononm}\frac{m}{m+\alpha(1+a_0)}\geq 1-\delta_2.\end{equation}

Case 1. $|h|^{2(1-\delta_2)}<t$, we choose $\epsilon=\frac{1}{2}|h|^{2(1-\delta_2)}$ and get (since $\delta_2m\geq \alpha(1+a_0)(1-\delta_2)$),

\begin{equation}\begin{aligned}\left|\mathbb{E}[\Delta_h^m\phi(X_t)]\right|\leq& \ell(t)(|h|^{\alpha(1+a_0)(1-\delta_2)}+|h|^{m\delta_2})\|\phi\|_{\mathscr{C}_b^\alpha}\\\leq&\ell(t)|h|^{\alpha(1+a_0)(1-\delta_2)}\|\phi\|_{\mathscr{C}_b^\alpha}.\end{aligned}\end{equation}

Case 2. $t\leq |h|^{2(1-\delta_2)},$ we choose $\epsilon=\frac{t}{2}$ and get (since $\delta_2<1$, we have that $m\geq \alpha(1+a_0)(1-\delta_2)$),
\footnote{The following computation is different from that in \cite{romito2018simple} as we have a different power of $1\wedge t$. We cannot use the computation given therein.}
\begin{equation}
    \begin{aligned}
    |\mathbb{E}[\Delta_h^m\phi(X_t)]|\leq& \ell(t)(t^{\frac{\alpha}{2}(1+a_0)}+\left(\frac{|h|}{\sqrt{t}}\right)^m)\|\phi\|_{\mathscr{C}_b^\alpha}\\\leq&\ell(t)(1\wedge t)^{-\frac{m}{2}}|h|^{\alpha(1+a_0)(1-\delta_2)}\|\phi\|_{\mathscr{C}_b^\alpha}.
    \end{aligned}
\end{equation}

In both cases, we have 
\begin{equation}\label{bothcases}
\left|\mathbb{E}[\Delta_h^m\phi(X_t)]\right|\leq \ell(t)(1\wedge t)^{-\frac{m}{2}}|h|^{\alpha(1+a_0)(1-\delta_2)}\|\phi\|_{\mathscr{C}_b^\alpha}.\end{equation}

Now we choose $\delta_2=\frac{\delta_1}{2\alpha(1+a_0)},$ we have $\delta_2\in(0,\frac{1}{2})$. Then $\alpha(1+a_0)(1-\delta_2)-\alpha=\alpha+\frac{\delta_1}{2}>\alpha,$ and as we have checked $m>\alpha(1+a_0)(1-\delta_2)$, we conclude from Lemma \ref{lemma4.1} that $p_x(t)\in B_{1,\infty}^a$. We also have the following qualitative estimate: for any $a<a_0$ we can find a decreasing function $\ell_a:[0,T]\to\mathbb{R}^+$ such that 
$$\sup_{s\geq t_0}\|p_x(s)\|_{B_{1,\infty}^a}\leq \ell_a(t_0)<\infty,\quad \text{ for each } t_0\in(0,T].$$
\end{proof}

In the proof of Corollary \ref{finalcorollary}, there is an additional assumption $\epsilon\leq h_t t<t$ for some $h_t\in(0,1).$ This is formally similar to the requirement $\epsilon\leq\frac{t}{2}$ in Proposition \ref{technicalproposition}. In light of this, in Case 1 we choose $\epsilon=(h_t\wedge \frac{1}{2})|h|^{2(1-\delta_2)}$, and in Case 2 we choose $\epsilon=(h_t\wedge \frac{1}{2})t.$  The overall difference we make in equation \eqref{bothcases} is a factor that only depends on $t$, while the power of $|h|$ is unchanged. Thus we still have the same conclusion.

To see why the estimates uniform over $s\geq t_0>0$ hold in Corollary \ref{finalcorollary}, we only need to note that the condition \eqref{secondsecondsecond} implies we can choose $h_t$ bounded away from $0$ uniformly over $t\in(t_0,T)$, for any $t_0>0$.

\section{Moments of the inverse Malliavin matrix} \label{momentappendix}

We review the proof of \cite{banos2016malliavin}, Proposition 4.4
and pay special attention to dependence of the estimate on $t$. 

Consider $X_t^n$ solution to $dX_t^n=a_n(t,X_t^n)dt+dW_t$ with drift $a_n\in \mathcal{C}^\infty$ with $a_n\to a$ a.e. such that $\sup_{n\geq 1}\|a_n'\|_\infty<M<\infty.$ 

Fix some $\delta\in(0,\frac{1}{2M})$ and consider $t>\delta$, using the series expansion  \eqref{theseriesexpansion} in the last inequality, we have
$$\|\mathscr{D}X_t^n\|^2_{L^2(\Omega;\mathbb{R}^{d\times d})}:=\int_0^t \|\mathscr{D}_sX_t^n\|_\infty^2ds\geq\int_{t-\delta}^t \|\mathscr{D}_sX_t^n\|_\infty^2ds\geq \frac{\delta}{2}-I_n(t,\delta)$$
where 
$$I_n(t,\delta):=\int_{t-\delta}^t \left\|\sum_{m\geq 1}\int_{s<u_1<\cdot<u_m<t}a_n'(u_1,X_{u_1}^n)\cdots a_n'(u_m,X_{u_m}^n)du_1\cdots du_m\right\|_\infty^2ds.$$

Since $a_n'$ are uniformly bounded by $M$ and $\delta\in(0,\frac{1}{2M})$, we can find a universal constant $C$ (not depending on $t$ and $\delta$ so long as $\delta\in(0,\frac{1}{2M})$) such that 
$$\sup_{n\geq 1}\mathbb{E}[|I_n(t,\delta)|^p]\leq C\delta^p,\quad p\geq 1.$$ Moreover, by choosing $\delta>0$ sufficiently small, an elementary computation shows that $\|\mathscr{D}X_t^n\|>0$ a.s. for each $t>0$.

Then by Chebyshev's inequality, for any $t>\delta$ and any $\epsilon<\frac{\delta}{2}$, the estimate
\begin{equation}
    \label{probestimates}
\begin{aligned}
\mathbb{P}\left(\|\mathscr{D}_\cdot X_t^n\|^2_{L^2(\Omega;\mathbb{R}^{d\times d})}\leq\epsilon\right)&\leq \mathbb{P}(I_n(t,\delta)\geq \frac{\delta}{2}-\epsilon)\\&\leq(\frac{\delta}{2}-\epsilon)^{-p}E[|I_n(t,\delta)|^p]\leq C\delta^p(\frac{\delta}{2}-\epsilon)^{-p}\end{aligned}\end{equation}
holds for any $p\geq 1$ and $\epsilon<\frac{\delta}{2}.$ We stress again that $C$ does not depend on $\delta\in(0,\frac{1}{2M}).$

The following inequality (see \cite{banos2016malliavin}, Lemma 4.3) is quite useful: for any random variable $Z:\Omega\to\mathbb{R}^d$ such that $\|Z\|>0$ a.s., for any $\eta_0>0$ and any $p\geq 1$, we have 
\begin{equation}
\label{estimatemoments}\mathbb{E}[\|Z\|^{-2p}]\leq \eta_0+\int_{\eta_0}^\infty \mathbb{P}(\|Z\|^{-2p}>\eta)d\eta= \eta_0+ p\int_0^{\eta_0^{-1/p}}\epsilon^{-(p+1)}\mathbb{P}(\|Z\|^2<\epsilon)d\epsilon,\end{equation}
where in the last equality we use the change of variable $\eta=\epsilon^{-p}.$

Now consider a test function (as in \cite{banos2016malliavin})
$$\delta(\epsilon):=\left|\frac{2\epsilon^{\frac{1}{2p}+2}}{\epsilon^{\frac{1}{2p}+1}-2}\right|.$$
It satisfies $\lim_{\epsilon\to 0+}\delta(\epsilon)=0$. Moreover, 
\begin{equation}\label{finiteness}\int_0^1 \epsilon^{-(p+1)}\left(\frac{\delta(\epsilon)}{2}-\epsilon\right)^{-p}\delta(\epsilon)^pd\epsilon<\infty\end{equation}
holds for any $p\geq 1$.

Now we conclude: for any $t_0>0$, we can find an upper bound of $\mathbb{E}[\|\mathscr{D}_\cdot X_t^n\|^{-p}]$ uniformly over all $t\geq t_0$, via the following procedure: (i) fix some $\delta_0\in(0,t_0\wedge \frac{1}{2M})$ and find an $\epsilon_0>0$ sufficiently small such that $\delta(\epsilon_0)\leq \delta_0;$ (ii) choose $\eta_0=\epsilon_0^{-p}$ in equation \eqref{estimatemoments}, set $Z=\mathscr{D}_\cdot X_t^n$, and plug in the bound \eqref{probestimates}; (iii) conclude via the finiteness estimate \eqref{finiteness}.

One may start the SDE $dX_t=a(t,X_t)dt+dW_t$ from any other point $x'\in\mathbb{R}^n$, then compute its Malliavin matrix and estimate moments of its inverse. Going through the same argument, one sees our quantitative estimate $\mathbb{E}[\|\mathscr{D}_\cdot X_t\|^{-p}]$ does not depend on the particular choice of initial value, i.e. can be uniform over all starting points $x'\in\mathbb{R}^n.$

\section{SDEs with non-smooth densities}\label{appendixB}

As promised in the introduction, we give an example of a one-dimensional SDE with additive Brownian noise, whose density function can be computed explicitly. The density function is Lipschitz, but no more regular than that. 
Consider the one-dimensional stochastic differential equation
$$dX_t=-\operatorname{sign}(X_t)dt+dW_t$$
with initial value $X_0=0$. The sign function is defined as $$\operatorname{sign}(x)=\begin{cases} 1\quad x\geq 0\\-1\quad x<0.
\end{cases}$$
Denote by $p_0(t)$ its density function at time $t$, where $0$ stands for the initial value. This SDE is closely related to the Brownian motion local times, via the following Tanaka's formula (see Exercise 5.3.12 of \cite{karatzas2012brownian}): for $W$ a Brownian motion, we have 
$$|W_t|=\int_0^t \operatorname{sign}(W_s)dW_s+2L_t^W(0),$$
where $L_t^W(0)$ is the local time at the origin for $W$. Therefore by Girsanov transform, let $\mathbb{Q}$ denote the law of the SDE, assuming $X$ has the law of a Brownian motion under $\mathbb{E}^0$, we have
$$\mathbb{Q}[X_t\in\Gamma]=e^{-t/2}\mathbb{E}^0[1_{X_t\in\Gamma}\exp(-|X_t|+2L(0))],\quad \text{ for any  measurable } \Gamma\subset\mathbb{R}.$$

Using explicit formulas for Brownian motion local times in dimension 1, the density of the SDE at $y\in\mathbb{R}$, which we denote by $p_0(t)(y)$, has the following explicit formula: (see for example Exercise 6.3.5 of \cite{karatzas2012brownian}): 
$$p_0(t)(y)=\frac{1}{\sqrt{2\pi t}}\left[\exp(\frac{(|y|+t)^2}{2t})+e^{-2|y|}\int_{|y|}^\infty \exp(-\frac{(v-t)^2}{2t})dv
\right].$$
This density function is Lipschitz but has no more regularity at $y=0$.

One may also consider the same SDE with initial value $x\in\mathbb{R}^d$, and compute the corresponding density function $p_x(t)(y)$. This is done in \cite{banos2014optimal}, Lemma 3.2:
$$p_x(t)(y)=\frac{1}{\sqrt{2\pi t}}e^{-\frac{(\operatorname{sign}(x)(x-y)-t)^2}{2t}}(1-e^{-\frac{2xy}{t}})1_{\operatorname{sign}(xy)\geq 0}+\int_0^t p_0(t-s)(y)\rho_{\tau_0^x}(s)ds,$$
where $\rho_{\tau_0^x}(s)$ is the first hitting time of $X_s$ to the origin, and thus has a smooth density. This time it is also clear that the density is Lipschitz continuous at $y=0$, but has no more regularity than that.

\section*{Acknowledgements}
The author thanks his supervisor, James Norris, for helpful suggestions on the paper.

\printbibliography

\end{document}